\documentclass[10pt]{article}
\usepackage{latexsym,amsmath,amssymb,graphics,amscd, empheq, mathtools}
\usepackage{enumitem, multirow}
\mathtoolsset{showonlyrefs}
\usepackage{xcolor}
\usepackage{stackengine}
\usepackage[breaklinks,colorlinks,
  citecolor=blue,
  urlcolor=blue]{hyperref}

\DeclareFontFamily{OT1}{pzc}{}
\DeclareFontShape{OT1}{pzc}{m}{it}{<-> s * [1.10] pzcmi7t}{}
\DeclareMathAlphabet{\mathpzc}{OT1}{pzc}{m}{it}

\addtocounter{footnote}{1}
\makeatletter
\@addtoreset{table}{bsection}
\def\thetable{\thesection.\@arabic\c@table}
\def\fps@table{h, t}
\@addtoreset{equation}{section}

\makeatother

\newtheorem{theorem}{Theorem}[section]
\newtheorem{definition}[theorem]{Definition}
\newtheorem{lemma}[theorem]{Lemma}
\newtheorem{remark}[theorem]{Remark}
\newtheorem{proposition}[theorem]{Proposition}

\newtheorem{corollary}[theorem]{Corollary}

\providecommand{\abs}[1]{\lvert#1\rvert}
\providecommand{\norm}[1]{\lVert#1\rVert}

\newcommand{\vertiii}[1]{{\left\vert\kern-0.25ex\left\vert\kern-0.25ex\left\vert #1 
    \right\vert\kern-0.25ex\right\vert\kern-0.25ex\right\vert}}
\newsavebox{\savepar}

\makeatletter

\usepackage{scalerel,stackengine}
\stackMath
\newcommand\reallywidehat[1]{%
\savestack{\tmpbox}{\stretchto{%
  \scaleto{%
    \scalerel*[\widthof{\ensuremath{#1}}]{\kern-.6pt\bigwedge\kern-.6pt}%
    {\rule[-\textheight/2]{1ex}{\textheight}}
  }{\textheight}%
}{0.5ex}}%
\stackon[1pt]{#1}{\tmpbox}%
}

\begin{document}

\title{\textbf{Forecasting causal dynamics with universal reservoirs}}
\author{Lyudmila Grigoryeva$^{1, 2}$, James Louw$^{3}$, and Juan-Pablo Ortega$^{3,4}$}
\date{}
\maketitle

\begin{abstract}
An iterated multistep forecasting scheme based on recurrent neural networks (RNN) is proposed for the time series generated by causal chains with infinite memory. This forecasting strategy contains, as a particular case, the iterative prediction strategies for dynamical systems that are customary in reservoir computing. Explicit error bounds are obtained as a function of the forecasting horizon, functional and dynamical features of the specific RNN used, and the approximation error committed by it. In particular, the growth rate of the error is shown to be exponential and controlled by the top Lyapunov exponent of the proxy system. The framework in the paper circumvents difficult-to-verify embedding hypotheses that appear in previous references in the literature and applies to new situations like the finite-dimensional observations of functional differential equations or the deterministic parts of stochastic processes to which standard embedding techniques do not necessarily apply.
\end{abstract}

\bigskip

\textbf{Key Words:} forecasting, dynamical system, ergodicity, nonuniform hyperbolicity, Oseledets' theorem, recurrent neural network, reservoir computing, state-space system, echo state network, ESN, generalized synchronization, causal chain with infinite memory, universality, universal approximation.

\makeatletter
\addtocounter{footnote}{1} \footnotetext{%
Mathematics and Statistics Division. University of St.~Gallen. Rosenbergstrasse 22,
St.~Gallen CH-9000, Switzerland.{\texttt{Lyudmila.Grigoryeva@unisg.ch}}}
\makeatother
\makeatletter
\addtocounter{footnote}{1} \footnotetext{%
Honorary Associate Professor at the Department of Statistics, University of Warwick, Coventry CV4 7AL, UK. {\texttt{Lyudmila.Grigoryeva@warwick.ac.uk} }}
\makeatother
\makeatletter
\addtocounter{footnote}{1} \footnotetext{%
Division of Mathematical Sciences,
School of Physical and Mathematical Sciences,
Nanyang Technological University,
21 Nanyang Link,
Singapore 637371.
{\texttt{Louw0002@e.ntu.edu.sg, 
Juan-Pablo.Ortega@ntu.edu.sg}}}
\makeatother
\makeatletter
\addtocounter{footnote}{1} \footnotetext{%
Corresponding author.}
\makeatother

\medskip

\medskip

\medskip

\section{Introduction}

{\it Recurrent neural networks} (RNN) and in particular {\it reservoir computing} (RC)  \cite{lukosevicius, tanaka:review} have gained much visibility in recent years due to their remarkable success in the learning of the chaotic attractors of complex nonlinear infinite dimensional dynamical systems \cite{Jaeger04, pathak:chaos, Pathak:PRL, Ott2018, pecora2020dependence, hart:ESNs, RC18, arcomano2022hybrid} and in a great variety of empirical classification and forecasting applications \cite{Wyffels2010, Buteneers2013, GHLO2014}). Even though much research has been conducted in recent years that theoretically explains this good performance, most of the work has been concentrated on the estimation \cite{RC10} and approximation  \cite{RC6, RC7, RC8, RC20, RC13, RC12}  properties of input/output systems but not in the way in which reservoir computing has been traditionally used as an RNN in order to forecast the observations of a dynamical system.

This paper aims to quantify how the universal approximation properties of a family of recurrent networks in the sense of input/output systems impact their approximation properties in the realm of the modeling and forecasting of the observations of autonomous dynamical systems and, more generally, of {\it causal chains with infinite memory} (CCIM) \cite{Dedecker2007a, doukhan2008weakly, alquier:wintenberger}. More specifically, we shall quantify the forecasting performance of an RNN that has been trained to mimic the observations of a dynamical system or a causal chain to a certain accuracy as an input/output system (where the output is a time-shifted version of the system observations in the input). This result will be obtained with minimal requirements on the state-space map. The upshot of our main result (Theorem~\ref{thm:total_state_error_bound}) is that the short-term forecasting (or path-continuation) error committed by an RNN that has been trained to uniformly approximate the dynamics of the observations up to a certain accuracy $\varepsilon>0$ is bounded above by the product of $\varepsilon$ with an exponential factor controlled by the top Lyapunov exponent of the autonomous dynamical system associated with the RNN that has been trained using the observations.

The two most popular forecasting schemes for causal dynamics are the so-called {\it iterated multistep} and {\it direct} forecasting strategies. There is a consensus on the fact that direct forecasting methods are superior in robustness with respect to model misspecification to those based on iterated multistep schemes such as those with RNNs \cite{Marcellino2006, Ghysels2009a, mcelroy2015direct}. This statement can be rigorously proved in many contexts, especially when estimation errors can be neglected. In our work, we focus on the iterated multistep approach since we have an interest not only in minimizing short-term forecasting errors but ultimately in  {\it using the trained system as a proxy for the underlying data-generating process}.

The causal chains that we forecast in this paper contain, as particular cases, two important dynamical elements. First, given a dynamical system $\phi \in \operatorname{Diff} ^1(M) $ on a $q$-dimensional manifold $M$ and $\omega \in C ^2(M, \mathbb{R}^d)$ a map that provides lower dimensional observations of it, Takens' Embedding Theorem \cite{takensembedding} proves that time series of $\omega  $-observations of $\phi$ form generically (that is, for an open and dense set of observation maps $\omega$ in the $C ^2 $ topology) a causal chain (actually of finite memory). Second, causal chains also appear as the deterministic parts of many stochastic processes after using the Wold \cite[Theorem 5.7.1]{BrocDavisYellowBook} or the Doob \cite{foellmer:schied:book} decompositions.

The observations-based modeling and forecasting in the dynamical systems context is a venerable topic that has been traditionally tackled using attractor reconstruction techniques based on the Takens' Embedding Theorem \cite{takensembedding, stark1999regularity, huke:2006, huke:muldoon:2015} and embedology \cite{sauer1991embedology}. More recently, new links have been established between traditional RNN-based forecasting schemes and dynamical systems learnability with classical notions like attractor embeddings and generalized synchronizations. The reader can check with \cite{hart:ESNs, lu:bassett:2020, manjunath:prsl, Verzelli2020, Verzelli2020b, allen:tikhonov, RC18, RC21, manjunath2021universal, manjunath2022embedding} and references therein for recent work in this direction. 

A particularly interesting recent work in the context of dynamical systems is \cite{berry2022learning}, which introduces a very general embedding framework that contains, as a particular case, most of the schemes in the references that we cited above. Moreover, it provides error estimates for the iterative multistep and direct forecasting strategies in that general framework, which are somewhat related to the results that we present in this paper. In that respect, we conclude this introduction by emphasizing some features of our approach that make it useful and, at the same time, differentiate it from the references above:

\begin{itemize}
\item We focus exclusively on the modeling and forecasting of the {\it observations} of a given dynamical system and {\it not of the full dynamical system itself}. This choice allows us to circumvent the need for embedding and general synchronization properties of the state-space system used for the learning and the dynamical system to be forecasted. This fact is particularly important because {\it it is difficult to establish that a given state-space system driven by prescribed observations of a dynamical system has embedding properties}, a hypothesis that appears in most references like \cite{berry2022learning}. The embedding property has been proved for certain linear Takens-like systems \cite{RC21} but remains to be shown for many standard universal RC and RNN paradigms like echo state networks \cite{Jaeger04} or state-affine systems \cite{RC6}, even though the synchronization maps that they induce have been proved to be typically differentiable \cite{RC18}.
\item Notwithstanding what has been said above, when given access to the full states of a dynamical system, our work also provides bounds for forecasting the dynamical system itself via two strategies, namely (i) seeing the identity map on the dynamical system as an observation map, and (ii) the causal embedding learning scheme introduced in \cite{manjunath2021universal}. The latter scheme uses a state map to learn a dynamical system from inputs consisting of the dynamical system's left-infinite trajectories. When the reservoir map satisfies a condition called {\it  state-input invertibility} (see \cite{manjunath2021universal} or Section~\ref{The causal embedding forecasting scheme}  later on for details), the scheme is a particular case of our setup.

\item Unlike the bounds for the iterated multistep scheme in \cite{berry2022learning}, the ones that we propose in this work are explicit (can be derived from constants introduced in the proofs) and are valid within a finite time horizon and not just asymptotically.
\item Our treatment covers not only observations of dynamical systems but also extends to general causal chains with infinite memory. We can hence handle forecasting problems that go beyond the situations usually treated with embedding techniques. These could be the solutions of infinite dimensional or functional differential equations (see example in Section~\ref{Numerical illustration}) or the deterministic parts of stochastic processes.

\item The rate at which the error grows with the forecasting horizon in the bounds is provided, and the dependence of the error bounds on the approximation error and analytic and dynamical features of the RNN used for forecasting is explicit. The dependence on the approximation error makes it possible to use bounds developed in previous works (see \cite{RC12, RC24}) that estimate it as a function of relevant architecture parameters in many popular reservoir families.
\end{itemize}

\paragraph{Outline of the paper.}

In Section~\ref{sec:preliminaries} we give an introduction to state-space systems and causal chains with infinite memory, as well as a brief overview of the ergodic theorems for dynamical systems needed in the proofs which follow later in the paper. Section~\ref{sec:Forecasting causal chains with recurrent neural networks} develops an analysis of the error in the iterated multistep forecasting scheme, from which we derive a bound for a simplified problem in Theorem~\ref{thm:linearized_seq_bound}, before arriving at our main result, Theorem~\ref{thm:total_state_error_bound}. Thereafter, we demonstrate how our results may be applied to the causal embedding scheme introduced in \cite{manjunath2021universal}. We conclude with Section~\ref{sec:numerial illustration}, in which we illustrate the bound \eqref{ineq:bound_deltay_final} for three different systems.

\paragraph{Codes and data.} All the codes and data necessary to reproduce the result in the paper are available at this \href{https://github.com/Learning-of-Dynamic-Processes/Forecasting-Causal-Chains-with-Universal-Reservoirs.git}{GitHub repository}.

\section{Preliminaries}\label{sec:preliminaries}

This section describes the dynamic processes that are the object of our research and the forecasting scheme we will use for them. Additionally, we spell out some classical results on ergodic theory that will be central in our developments.

\paragraph{The data generating process (DGP) and its approximants.} 

Our results apply to the discrete-time evolution of dynamical variables in a $d$-dimensional Euclidean space whose behavior is determined by a {\it causal chain with infinite memory} (CCIM). More specifically, we are interested in the dynamical properties of infinite sequences ${\bf y}= \left({\bf y}_t\right)_{t \in \mathbb{Z}}$, where ${\bf y}_t\in \mathbb{R}^d $, for all $t \in \mathbb{Z} $, for which there exists a functional $H: (\mathbb{R}^d)^{\mathbb{Z}_{-}} \longrightarrow \mathbb{R}^d $ such that 
\begin{equation*}
{\bf y} _t=H \left(\ldots,  {\bf y} _{t-2}, {\bf y} _{t-1}\right)=H (\underline{{\bf y} _{t-1}}), \quad \mbox{for all $t \in \mathbb{Z} $,} 
\end{equation*}
and where the symbol $\underline{{\bf y} _{t-1}} \in  ({\mathbb{R}}^d)^{\mathbb{Z}_{-}}$ denotes the semi-infinite sequence $\underline{{\bf y} _{t-1}}=\left(\ldots,  {\bf y} _{t-2}, {\bf y} _{t-1}\right) $. These dynamics are the deterministic counterparts of various structures with a widespread presence in time series analysis. See for instance \cite{Dedecker2007a, doukhan2008weakly, alquier:wintenberger} and references therein. Causal chains also have strong connections to {\it symbolic dynamics} \cite{Das2024}. The (left) shift map $\sigma \colon \mathcal{Y}^{\mathbb{Z}} \longrightarrow \mathcal{Y}^{\mathbb{Z}}, \, \sigma({\bf y})_t = {\bf y}_{t+1}$ is related to the functional $H$ associated to the CCIM by

\begin{align*}
    H \circ \pi_{\mathbb{Z}_-} \circ \sigma^t ({\bf y}) = \pi_0 \circ \sigma^t({\bf y}) = \pi_t({\bf y}), \text{ for all } t \in \mathbb{Z},
\end{align*}

where $\pi_{\mathbb{Z}_-}\colon \mathcal{Y}^{\mathbb{Z}} \longrightarrow \mathcal{Y}^{\mathbb{Z}_-}$ is the projection onto the left infinite history, and $\pi_t \colon \mathcal{Y}^{\mathbb{Z}} \longrightarrow\mathcal{Y}$ the projection onto the $t$-th coordinate. Alternatively, we may write

\begin{align*}
    (\dots, {\bf y}_{t-2}, {\bf y}_{t-1}, H(\underline{{\bf y}_{t-1}})) = \pi_{\mathbb{Z}_-} \circ \sigma^{t+1} ({\bf y}) \text{ for all } t \in \mathbb{Z}.
\end{align*}

Thus the family of sequences that form the causal chain are a shift invariant subspace of $\mathcal{Y}^{\mathbb{Z}}$. While in symbolic dynamics one typically considers shift spaces over a finite alphabet, we allow the space $\mathcal{Y}$ to be more general, possibly uncountable. Typically we take $\mathcal{Y} = \mathbb{R}$. In the field of dynamical systems, symbolic dynamics supplies a very tractable means of studying integral concepts such as topological entropy, and more generally, ergodic theory \cite{Katok_Hasselblatt_1995}. A number of key examples lending insight into the properties of such quantities come from symbolic dynamics. Certain classes of dynamical systems may be approximated by shift spaces through the use of Markov partitions. The application of these ideas to general causal chains would be an interesting direction of study.

An important aspect of causal chains, which poses problems in carrying out numerical realizations of such chains, is that they depend on an infinite history of values, whereas a computer algorithm can only take as input a finite length sequence. This raises the question of the {\it computability} of causal chains. While we make some comments below on the approximation of causal chains by functions with finite length inputs, readers may look at \cite{Lind_Marcus_1995,kitchens2012symbolic} for an introduction to computability theory. A further interesting reference discussing the interaction between symbolic dynamics and computability theory may be found in \cite{Burr_2022}.

The CCIM determined by the functional $H: ({\mathbb{R}}^d)^{\mathbb{Z}_{-}} \longrightarrow {\mathbb{R}}^d $ has a canonical {\it causal and time-invariant filter }$U _H: ({\mathbb{R}}^d)^{\mathbb{Z}} \longrightarrow ({\mathbb{R}}^d)^{\mathbb{Z}} $ associated determined by the equality
\begin{equation}
\label{causal and tinv filter}
U _H( {\bf y})_t:=H (\underline{{\bf y} _{t-1}}), \quad \mbox{for all ${\bf y} \in ({\mathbb{R}}^d)^{\mathbb{Z}}$ and  $t \in \mathbb{Z} $.} 
\end{equation}

A classical problem in control and systems theory is the {\it realization} of filters like \eqref{causal and tinv filter} as the unique solution (when it exists) of a state-space transformation
\begin{empheq}[left={\empheqlbrace}]{align}
\mathbf{x} _t &=F(\mathbf{x}_{t-1}, {\bf y} _{t-1}),\label{eq:RCSystem1}\\
{\bf y}_t &= h (\mathbf{x} _t), \label{eq:RCSystem2}
\end{empheq}
where $F: \mathbb{R}^L \times \mathbb{R} ^d \longrightarrow \mathbb{R}^L $ and $h : \mathbb{R}^L \longrightarrow \mathbb{R}^d$ are a {\it state} and a  {\it readout} map, respectively, $\mathbb{R}^L  $ is, in this case, the state space and, typically, $N$ is much
bigger than $d$. Some classical references in connection to this problem are \cite{kalman1959general, kalman1959unified, kalman:original, Kalman1962, baum1966statistical, Kalman2010}  or, more recently, \cite{kalman:mit, Sontag1979, DangVanMien1984, Narendra1990, Matthews1994, lindquist:picci, RC16, RCSP1}. Even though the realization problem does not always have a solution, it has been shown that if the domain of the functional $H $ is restricted to a space of uniformly bounded sequences of the type $K _M:= \left\{{\bf y}\in (\mathbb{R}^d)^{\mathbb{Z}_-}\mid \left\|{\bf y} _t\right\|\leq M, \, t \in \mathbb{Z}_{-}\right\} $ for some $M >0 $ and, additionally, it has a continuity condition called the {\it fading memory property} \cite{Boyd1985, RC6}, then various families of state-space transformations of the type \eqref{eq:RCSystem1}-\eqref{eq:RCSystem2} can be used to approximate it uniformly. This is the case for {\it linear systems with polynomial or neural network readouts} \cite{Boyd1985, RC6, RC8}, {\it state-affine systems with linear readouts} \cite{RC6}, or {\it echo state networks} \cite{RC7, RC8, RC12, RC20}. In other words, any of these families exhibit {\it universality properties} in the fading memory category. This allows us to conclude that given a fading memory functional $H: K _M \longrightarrow \mathbb{R}^d $ and $\varepsilon>0  $, there exists a state-space system $(F , h) $ in that family to which a functional $H^F_h:K _M \longrightarrow \mathbb{R}^d $ can be associated and such that
\begin{equation*}
\left\|H-H^F_h\right\| _{\infty}=\sup_{{\bf y}\in K _M} \left\{\left\|H({\bf y})-H^F_h({\bf y})\right\|\right\}< \varepsilon.
\end{equation*}
For future reference, a condition that ensures that a state-space system $(F , h) $ has a functional $H^F_h:K _M \longrightarrow \mathbb{R}^d $ associated is the so-called {\it echo state property} (ESP) \cite{jaeger2001}: the state maps $F $ has the ESP with respect to inputs in $K _M $ when for any ${\bf y} \in K _M $ there exists a unique $\mathbf{x} \in (\mathbb{R}^L)^{\mathbb{Z}_{-}} $ such that  \eqref{eq:RCSystem1} holds for any $t \in \mathbb{Z}_{-} $. In such a situation, $F$ has a filter $U ^F: K _M \longrightarrow (\mathbb{R}^L)^{\mathbb{Z}_{-}} $ associated that is uniquely determined by the relations $U ^F({\bf y})_t=F \left(U ^F({\bf y})_{t-1}, {\bf y}_{t-1}\right) $. The filter $U^F$ induces the functional $H^F \colon K_M \longrightarrow {\mathbb{R}}^N$ given by $H^F({\bf y}) = U^F({\bf y})_0$. The functional $H^F_h:K _M \longrightarrow \mathbb{R}^d $ corresponding to the state-space system is given by $H^F_h({\bf y})=h \circ H^F({\bf y})$. The ESP has important dynamical implications \cite{RC9, manjunath:prsl, manjunath2022embedding}. In practice, the state space is often taken to be compact. In this case, if $F$ is a continuous map and a contraction on its first entry, the ESP is ensured. Indeed, for such a state map, it can be shown that with inputs from $\overline{B({\bf 0}, M)} \subset \mathbb{R}^d$ there is a compact subset $X \subset \mathbb{R}^L$ such that $F(X \times \overline{B({\bf 0}, M)}) \subset X$. Stochastic generalizations of some of these statements can be found in \cite{RC27, RC28}.

\paragraph{Causal chains with infinite memory and dynamical systems observations.} 

A particularly important family of CCIMs can be generated using the observations of some dynamical systems. Indeed, let $M$ be a manifold and let $\phi:M \longrightarrow M $ be a continuous map that determines a dynamical system on it via the evolution recursions
\begin{equation*}
m _{t+1}= \phi(m _t),
\end{equation*}
using a given point $m _0 \in M $ as initial condition. The data available for the forecasting of a dynamical system is given in the form of observations that are usually encoded as the values of a continuous map $\omega:M \longrightarrow {\mathbb{R}}^d $, that is, they are given by the sequences ${\bf y}_0:= \omega(m _0),\, {\bf y}_1:= \omega(m _1), \ldots, {\bf y}_t:= \omega(m _t)$. The problem that we tackle in this work is the forecasting or continuation of those sequences using state-space transformations of the type  \eqref{eq:RCSystem1}-\eqref{eq:RCSystem2} that have been learned/trained using exclusively temporal traces of $\omega $-observations.

We shall now show that the observations of a large class of dynamical systems are CCIMs. More specifically, we shall work in the context of {\it invertible dynamical systems} for which an {\it invertible generalized synchronization} map can be found. We briefly recall the meaning of these terms. First, we say that the dynamical system $\phi $ on $M$ is invertible when it belongs either to the set  ${\rm Hom}(M)$ of homeomorphisms (continuous invertible maps with continuous inverse) of a topological space $M$ or to the set of diffeomorphisms  ${\rm Diff} ^1(M) $ of a finite-dimensional manifold $M $. The invertibility of $\phi$  allows us, given an observation map $\omega: M \longrightarrow {\mathbb{R}}^d $, to define  the $(\phi, \omega)$-{\it delay map} $S _{(\phi, \omega)}:M \longrightarrow (\mathbb{R}^d)^{\mathbb{Z}} $ as $S _{(\phi, \omega)}(m):=\left(\omega(\phi ^t (m))\right)_{t \in \mathbb{Z}} $.  Second, given a state-space map $F: \mathbb{R} ^L\times \mathbb{R} ^d \longrightarrow   \mathbb{R} ^L$, consider the drive-response system associated to the $\omega $-observations of $\phi$ and determined by the recursions:
\begin{equation}
\label{drive-response system}
\mathbf{x} _t=F\left(\mathbf{x} _{t-1}, S _{(\phi, \omega)}(m)_{t-1}\right), \quad \mbox{$t \in \mathbb{Z},\, m \in M.$}
\end{equation}
We say that a {\it generalized synchronization} (GS)  \cite{rulkov1995generalized, pecora:synch, ott2002chaos, boccaletti:reports:2002, eroglu2017synchronisation}  occurs in this configuration when there exists a map $f_{(\phi, \omega, F)}:M \longrightarrow \mathbb{R}^L $ such that for any $\mathbf{x} _t  $, $t \in \mathbb{Z} $, as in~\eqref{drive-response system}, it holds that
\begin{equation}
\label{generalized synchronization condition}
 \mathbf{x} _t = f_{(\phi, \omega, F)} (\phi ^t(m)),
\end{equation}
that is, the time evolution of the dynamical system in phase space (not just its observations) drives the response in~\eqref{drive-response system}. We emphasize that the definition \eqref{generalized synchronization condition} presupposes that the recursions \eqref{drive-response system} have a (unique) solution, that is, that there exists a sequence $\mathbf{x} \in (\mathbb{R}^L)^{\mathbb{Z}} $ such that \eqref{drive-response system} holds true. When that existence property holds and, additionally, the solution sequence $\mathbf{x}  $ is unique, we say that $F$ has the $(\phi, \omega) $-{\it Echo State Property} (ESP). The existence, continuity, and differentiability of GSs have been established, for instance, in \cite{stark1999regularity, RC18} for a rich class of systems that exhibit the so-called fading memory property and that are generated by locally state-contracting maps $F$. In such a framework, the state equation \eqref{drive-response system} uniquely determines a filter $U^F: \left({\mathbb{R}}^d\right)^{\mathbb{Z}} \longrightarrow \left({\mathbb{R}}^L\right)^{\mathbb{Z}}$  that satisfies the recursion $U^F({\bf z}) _t=F\left(U^F({\bf z}) _{t-1}, {\bf z} _{t-1}\right)$, for all $t \in \mathbb{Z}$, in terms of which the GS $f_{(\phi, \omega, F)}:M \longrightarrow \mathbb{R}^L $ can be written as
\begin{equation}
\label{charact GS}
f_{(\phi, \omega, F)}(m)= \left(U^F \left(S _{(\phi, \omega)}(m)\right)\right)_0.
\end{equation}

The following result, introduced in \cite{RC18, RC26}, shows that the observations of a dynamical system in the presence of injective GS form a CCIM. More explicitly, we shall be assuming that the GS $f_{(\phi, \omega, F)}:M \longrightarrow \mathbb{R}^L $ is an {\it embedding} which means that $f_{(\phi, \omega, F)}$  is an injective immersion ($C ^1  $ map with injective tangent map) and, additionally, the manifold topology in $f_{(\phi, \omega, F)}(M) $ induced by $f_{(\phi, \omega, F)}$ coincides with the relative topology inherited from the standard topology in $\mathbb{R}^L $. 

\begin{proposition}
\label{GS and conjugacies}
Let $\phi:M \longrightarrow M $ be an invertible dynamical system on the manifold $M$ and let $\omega:M \longrightarrow \mathbb{R}^d  $ be an observation map. Let $F: \mathbb{R} ^L\times \mathbb{R} ^d \longrightarrow   \mathbb{R} ^L $ be a state map that determines a generalized synchronization $f_{(\phi, \omega, F)}:M \longrightarrow \mathbb{R}^L $ when driven by the $\omega $-observations of $M$ as in~\eqref{drive-response system}. Suppose that $f_{(\phi, \omega, F)}:M \longrightarrow \mathbb{R}^L $ is an embedding. Then:
\begin{description}
\item [(i)] The set $S:=f_{(\phi, \omega, F)}(M) \subset \mathbb{R}^L  $ is an embedded submanifold of the reservoir space  $\mathbb{R}^L  $.
\item [(ii)] There exists a differentiable observation map $h: S\longrightarrow \mathbb{R}^d  $ that extracts the corresponding observations of the dynamical system out of the reservoir states. That is, with the notation introduced in~\eqref{drive-response system}:
\begin{equation}
\label{readout forecast}
h(\mathbf{x} _t)= \omega\left(\phi ^{t} (m)\right).
\end{equation}
\item [(iii)]  The maps $F$ and $h$ determine a differentiable dynamical system $\Phi \in C^1(S, S) $ given by 
\begin{equation*}
\label{conjugate dynamical system}
\Phi( {\bf s}):=F({\bf s}, h({\bf s})),
\end{equation*}
which is $C ^1 $-conjugate to $\phi\in {\rm Diff}^1(M) $ by $f_{(\phi, \omega, F)}$, that is,
\begin{equation*}
\label{c1 conjugation}
f_{(\phi, \omega, F)} \circ \phi = \Phi \circ f_{(\phi, \omega, F)}.
\end{equation*}
\item [(iv)] The sequences of $\omega $-observations ${\bf y}_t:= \omega(\phi^t(m))$, $t \in \mathbb{Z} $, for a fixed $m \in M $, are causal chains with infinite memory.
\end{description}
\end{proposition}

\noindent The proof of this result can be found in \cite{RC26} and is an elementary consequence of the fact that $f_{(\phi, \omega, F)}$ is an embedding (see, for instance, \cite{mta} for details). The map in {\bf (ii)} is $h:= \omega \circ f_{(\phi, \omega, F)} ^{-1}: f_{(\phi, \omega, F)}(M) \subset  \mathbb{R}^L \longrightarrow \mathbb{R} ^d $ (notice that it can be constructed since, by hypothesis, $f_{(\phi, \omega, F)}$ is invertible). Regarding {\bf (iv)}, this is a consequence of \eqref{readout forecast}, \eqref{generalized synchronization condition}, and \eqref{charact GS}:
\begin{equation*}
{\bf y}_t= \omega\left(\phi ^{t} (m)\right)=h(\mathbf{x} _{t})= h(f_{(\phi, \omega, F)} (\phi ^{t}(m)))= h \left(\left(U^F \left(S _{(\phi, \omega)}(\phi ^{t}(m))\right)\right)_0\right)=:H ^F_h( \underline{ {\bf y}_{t-1}}),
\end{equation*}
where the dependence of the newly defined map $H ^F_h $ exclusively on $\underline{ {\bf y}_{t-1}} $ is a consequence of the causality of the filter $U ^F $. 

\medskip

As we show in the next corollary, the last part of this proposition, together with Takens' Embedding Theorem (or, more generally, the main result in \cite{RC21}), implies that generic observations of a large class of dynamical systems form CCIMs.

\begin{corollary}\label{corr:dynamical_systems_and_causal_chains}
Let $M$ be a compact manifold of dimension $q \in \mathbb{N}$ and let $\phi \in {\rm Diff} ^2(M) $ be a twice-differentiable diffeomorphism that satisfies the following two properties:
\begin{description}
\item [(i)] $\phi $ has only finitely many periodic points with periods less than or equal to $2q$.
\item [(ii)] If $m \in M $  is any periodic point of $\phi $ with period $k<2q $, then the eigenvalues of the linear map $T _m \phi ^k: T _m M \longrightarrow T _mM  $ are distinct.
\end{description}
Then for any generic scalar observation function $\omega \in C ^2(M, \mathbb{R} )$, the sequences of $\omega $-observations are causal chains with infinite memory.
\end{corollary}

\noindent {\bf Proof.} Recall first that in the hypotheses of the corollary, Takens' Theorem \cite{takensembedding, huke:2006} guarantees that a $(2q+1)$-truncated version $S _{(\phi, \omega)} ^{2q+1}$ of the $(\phi, \omega)$-{delay map} given by  $S _{(\phi, \omega)}^{2q+1}(m):=\left(\omega (m), \omega(\phi ^{-1} (m)), \ldots, \omega(\phi ^{-2q} (m))\right)$ is a continuously differentiable embedding. This map is in turn, the GS corresponding to the linear state map $F(\mathbf{x}, z):=A \mathbf{x}+ \mathbf{C} z $, with $A$ the lower shift matrix in dimension $2q+1  $ and $\mathbf{C}= (1,0, \ldots,0) \in \mathbb{R}^{2q+1} $. The claim then follows from part {\bf (iv)} of Proposition~\ref{GS and conjugacies}. \quad $\blacksquare$

\paragraph{Ergodic theorems for dynamical systems.} An important tool in the study of dynamical systems is the concept of ergodicity, which allows us to translate statements about space averages to averages along trajectories of the system, thus relating global analytic properties of a map to the dynamics it induces. Of particular relevance to our work is the study of nonuniform hyperbolicity and the multiplicative ergodic theorem of Oseledets \cite{oseledets1968multiplicative,ruelle1979ergodic,barreira2007nonuniform}. It is this theorem which forms the cornerstone of the proof, yielding the exponential factor in the final bound \eqref{ineq:bound_deltay_final}. We briefly introduce the theory as applied to our context, keeping to the details and results necessary for our work. Our overview is based largely on \cite[Chapter 1 and Chapter 3]{barreira2007nonuniform}.

Let $M$ be an $q$-dimensional $C^1$ submanifold of $\mathbb{R}^{L}$ and let $\phi\in \operatorname{Diff}^1(M)$ be a dynamical system acting on $M$. A concept that allows us to quantify the rate of divergence of two points initially close together is that of the {\it Lyapunov exponent}. For a point $m \in M$ the {\it forward Lyapunov exponent} along a vector $v \in T_mM$ in the tangent space $T_mM$ of $M$ at the point $m$ is defined as
\begin{align}\label{LE_forward}
\lambda^+ (m, v) = \limsup_{t \to \infty} \frac{1}{t} \log \norm{D\phi^t(m)v},
\end{align}
where $D\phi^t(m):T _mM \longrightarrow T _{\phi^t(m)}M$ is the tangent map of $\phi^t $ at the point $m\in M$.
Similarly, we define the {\it backward Lyapunov exponent}
\begin{align}\label{LE_backward}
\lambda^- (m, v) = \limsup_{t \to -\infty} \frac{1}{|t|} \log \norm{D\phi^t(m)v}.    
\end{align}
It can be shown that, for any given point $m \in M$, the forward Lyapunov exponent $\lambda^+(m, \cdot)$ satisfies the following three properties:

\begin{description}
\item[(i)] $\lambda^+(m, \alpha v) = \lambda^+(v)$ for all $v \in T_mM$ and all $\alpha \in \mathbb{R} \setminus \{0\},$
\item[(ii)] $\lambda^+(m, v+w) \leq \max\{\lambda^+(m, v), \lambda^+(m, w)\}$ for all $v, w \in T_mM$, and 
\item[(iii)] $\lambda^+(0) = - \infty$, where we make use of the convention $\log 0 = - \infty.$
\end{description}
As a result, $\lambda^+(m, \cdot)$ takes on finitely many values $\lambda_1^+(m)>\cdots> \lambda_{l^+(m)}^+(m)$ for some $l^+(m) \leq q$. These are related to the differing rates of convergence/divergence along different tangent directions. The properties {\bf (i)}--{\bf (iii)} also allow us to construct a filtration of vector subspaces in the tangent space using the ordered Lyapunov exponents to single out linearly independent tangent directions up to a given Lyapunov value. The filtration is $\mathcal{V}^+(m)\colon \{ 0 \} \subsetneq V_{l^+(m)}^+(m) \subsetneq \cdots \subsetneq V_1^+(m)=T_mM$, where we define
\begin{align}\label{filtration_forward}
V_i^+(m) := \{ v \in T_mM\colon \lambda^+(m,v) \leq \lambda_i^+(m)\},
\end{align}
for $i=1,\dots,l^+(m)$, that is, $V_i^+(m)$ is the subspace of the tangent space consisting of all directions with asymptotic rate of divergence at most $\lambda_i^+(m)$. The number of tangent directions having a given rate of divergence is $k_i^+(m) = \dim V_i^+(m) - \dim V_{i+1}^+(m)$ for $i=1, \dots, l^+(m)$. The value $k_i^+(m)$ is referred to as the multiplicity of the Lyapunov exponent $\lambda_i^+(m).$

In like manner, for $m \in M$ the backward Lyapunov exponent $\lambda^-(m, \cdot)$ takes on finitely many values $\lambda_1^-(m)<\cdots< \lambda_{l^-(m)}^-(m)$ for some $l^-(m) \leq q$ and has a filtration $\mathcal{V}^-(m)\colon \{ 0 \} \subsetneq V_{1}^-(m) \subsetneq \cdots \subsetneq V_{l^-(m)}^-(m)=T_mM$, defined by
\begin{align*}
V_i^-(m) := \{ v \in T_mM\colon \lambda^-(m,v) \leq \lambda_i^-(m)\},
\end{align*}
for $i=1,\dots,l^-(m)$. In a similar fashion to the forward Lyapunov exponent, the multiplicity of the backward Lyapunov exponent $\lambda_i^-(m)$ is defined as $k_i^-(m) = \dim V_i^-(m) - \dim V_{i-1}^-(m)$ for $i=1, \dots, l^-(m).$

While the forward Lyapunov exponent captures the growth of tangent vectors in forward time, the backward exponent measures this growth in backward time. An obstacle in the definitions \eqref{LE_forward} and \eqref{LE_backward} is the $\limsup$, which makes the growth of individual terms in the sequences intractable. To deal with this, we relate the forward and backward Lyapunov exponents and their corresponding filtrations. This relation is expressed in the notion of {\it coherence.}
The filtrations $\mathcal{V}^+(m)$ and $\mathcal{V}^-(m)$ are said to be {\it coherent} if
\begin{description}
\item[(i)] $l^+(m) = l^-(m) = : l(m),$
\item[(ii)] there is a decomposition $T_mM = \bigoplus_{i=1}^{l(m)} E_i(m)$ such that
\begin{align}\label{Oseledets splitting}
V_i^+(m) = \bigoplus_{j=i}^{l(m)} E_j(m) \quad \text{and} \quad V_i^-(m)= \bigoplus_{j=1}^{i} E_j(m) \quad \text{for } i=1,\dots,l(m), \text{ and}
\end{align}
\item[(iii)] $\lim_{t \to \pm \infty} \frac{1}{t} \log \norm{D\phi^t(m)v} = \lambda_i^+(m) = - \lambda_i^-(m) =: \lambda_i(m)$ uniformly over $\{v \in E_i(m) \colon \norm{v} = 1\}.$
\end{description}
The subspaces $E_i(m)$ are called the {\it Oseledets subspaces} and the decomposition $T_mM = \bigoplus_{i=1}^{l(m)} E_i(m)$ is called the {\it Oseledets splitting}. The number of distinct Lyapunov exponents and the Oseledets subspaces are invariant along a trajectory of the dynamical system in the sense that $l(m) = l(\phi(m))$ and $D\phi(m) E_i(m) = E_i(\phi(m))$ for $i = 1, \dots, l(m)$. The Oseledets subspaces can be constructed by taking $E_i(m) = V_i^+(m) \cap V_i^- (m)$ for $i = 1, \dots, l(m)$. Expressing them in this way is instructive: the set $V_i^+(m)$ filters out those tangent directions which, in forward time, have greater Lyapunov exponent than $\lambda_i(m)$, while the set $V_i^-(m)$ filters out those that, in backward time, have greater Lyapunov exponent than $-\lambda_i(m)$. Thus the tangent directions left in $E_i(m)$ are in some sense {\it pure}, in that their Lyapunov exponent in forward and backward time are the same.

We focus on a certain {\it regular} class of points whose properties allow us to establish a link between the Lyapunov exponent and the rate of divergence of points in a nearby neighborhood using the theory of nonuniform hyperbolicity. In particular, for these points the $\limsup$s in~\eqref{LE_forward} and \eqref{LE_backward} become limits and they possess a splitting as in~\eqref{Oseledets splitting}, allowing us to decompose the tangent space into subspaces where the rate of expansion/contraction under the action of the dynamical system is uniform. A point $m \in M$ is said to be {\it forward regular} if $$\lim_{t \to \infty} \frac{1}{t} \log \abs{\det D \phi^t(m)} = \sum_{i=1}^{l^+(m)} k_i^+(m) \lambda_i^+(m) < \infty.$$ Likewise, $m$ is said to be {\it backward regular} if $$\lim_{t \to -\infty} \frac{1}{|t|} \log \abs{\det D \phi^t(m)} = \sum_{i=1}^{l^-(m)} k_i^-(m) \lambda_i^-(m) < \infty.$$

In both cases, the expression on the right is the sum of the distinct values of the Lyapunov exponent counted with multiplicity. Intuitively, this says that a small volume of space around the point $m$ grows at an exponential rate controlled by the sum of the Lyapunov exponents. Finally, we are ready to define the set of regular points: a point is said to be {\it Lyapunov-Perron regular}, or simply {\it regular}, if it is both forward and backward regular and the filtrations $\mathcal{V}^+(m)$ and $\mathcal{V}^-(m)$ are coherent. This brings us to {\it Oseledets' Multiplicative Ergodic Theorem} \cite[Theorem 3.4.3]{barreira2007nonuniform}.


\begin{theorem}[Oseledets' Multiplicative Ergodic Theorem]\label{thm:Oseledets_mult_erg_thm}
Let $\phi\in \operatorname{Diff}^1(M)$ be a dynamical system acting on the $q$-dimensional $C^1$ submanifold $M$ of $\mathbb{R}^{L}$. Let $\mu$ be a Borel probability measure on $M$, invariant with respect to $\phi$, and suppose $\int_M \operatorname{log}^+ \| D\phi \| d\mu<\infty$ and $\int_M \operatorname{log}^+ \| D\phi^{-1} \| d\mu < \infty$, where $\operatorname{log}^+ a = \max \{ \operatorname{log} a, 0\}$. Then the set $\Lambda \subset M$ of regular points is invariant under $\phi$ and has full $\mu$-measure.
\end{theorem}

While above we gave the standard definition for a regular point, we require a slightly stronger consequence of regularity found in \cite[Theorem 1.3.11]{barreira2007nonuniform}. The result involves a certain type of basis which is conveniently adapted to the associated filtration at the point.

\begin{definition}
Consider an $q$-dimensional vector space $V$ and a filtration $\{0\} = V_{l+1} \subsetneq V_{l} \subsetneq \cdots \subsetneq V_1 = V$. A basis $\{v_1, \dots, v_q\}$ of $V$ is said to be normal with respect to this filtration if for every $i = 1, \dots, l$ there is a basis for $V_i$ in the set $\{v_1, \dots, v_q\}.$
\end{definition}

The result we need is in fact an equivalent condition for forward regularity of a point, quantifying how a small volume at a point expands along trajectories. It allows us to bound how the angles between the Oseledets subspaces in the Oseledets splitting evolve with time. For a set $\{v_1, \dots, v_r\}$ of linearly independent vectors, we denote by $\operatorname{vol}(v_1, \dots, v_r)$ the volume of the $r$-parallelepiped formed by these vectors. The following result can be found in \cite[Theorem 1.3.11]{barreira2007nonuniform}.

\begin{theorem}
\label{thm:barr_pes_Oseledets_extension}
Let $M$ be an $q$-dimensional $C^1$ submanifold of $\mathbb{R}^{L}$ and let $\phi\in \operatorname{Diff}^1(M)$ be a dynamical system acting on $M$. Consider a point $m \in M$, with associated filtration $\mathcal{V}(m)\colon \{ 0 \} = V_{l(m)+1}(m) \subsetneq V_{l(m)}(m) \subsetneq \cdots \subsetneq V_1(m) = T_mM$, as defined above. Then $m$ is forward regular if and only if for any basis $\{v_1, \dots, v_q\}$ for $T_mM$, normal with respect to $\mathcal{V}(m)$, and for any $K\subset \{1, \dots, q\},$
\begin{align*}
\lim_{t \to \infty} \frac{1}{t} \log \operatorname{vol}(\{D\phi^t(m) v_i\}_{i \in K}) = \sum_{i \in K} \lambda(m,v_i).
\end{align*}
\end{theorem}




\section{Forecasting causal chains with recurrent neural networks}\label{sec:Forecasting causal chains with recurrent neural networks}

We now study the iterated multistep forecasting scheme for causal chains using recurrent neural networks and prove our main result. The main theorem provides error bounds for the forecasting error as a function of the forecasting horizon and the dynamical properties of the RNN used in the prediction exercise.

\subsection{The iterated multistep forecasting scheme} 

We formulate the forecasting problem in a simplified setup in which we are given a bi-infinite sequence ${\bf y}= \left({\bf y}_t\right)_{t \in \mathbb{Z}}$ which is a realization of a causal chain with infinite memory determined by a functional  $H: ({\mathbb{R}}^d)^{\mathbb{Z}_{-}} \longrightarrow {\mathbb{R}}^d $. This means that the entries of ${\bf y} $ satisfy ${\bf y} _t=H \left(\ldots,  {\bf y} _{t-2}, {\bf y} _{t-1}\right)=H \left(\underline{{\bf y} _{t-1}}\right)$, for all $t \in \mathbb{Z} $. The iterated multistep forecasting of ${\bf y} $ consists in producing predictions $\widetilde{ {\bf y}} _1, \widetilde{ {\bf y}} _2, \ldots, \widetilde{ {\bf y}} _t$ of the positive entries $ {\bf y}_1, {\bf y}_2 \ldots, {\bf y}_t  $, $t \in \mathbb{N} $ of ${\bf y} $ based on the knowledge of the semi-infinite sequence $\underline{{\bf y}_0} \in ({\mathbb{R}}^d)^{\mathbb{Z}_{-}} $. The RNN-based iterated multistep forecasting scheme consists of first finding a state-space system $(F,h) $ with associated functional $H^F_h:({\mathbb{R}}^d)^{\mathbb{Z}_{-}} \longrightarrow {\mathbb{R}}^d $ that uniformly $\varepsilon $-approximates $H$, that is,
\begin{equation}
\label{condition unif approx}
\left\|H-H^F_h\right\| _{\infty}=\sup_{{\bf y}\in ({\mathbb{R}}^d)^{\mathbb{Z}_{-}}} \left\{\left\|H({\bf y})-H^F_h({\bf y})\right\|\right\}< \varepsilon.
\end{equation}

Once such $(F,h)$ has been found, the predictions $\widetilde{ {\bf y}} _1, \widetilde{ {\bf y}} _2, \ldots, \widetilde{ {\bf y}} _t$ are iteratively produced in tandem with a sequence $\widetilde{ {\bf x}}_1, \widetilde{{\bf x}}_2, \ldots, \widetilde{ {\bf x}}_t$ in the state space according to the scheme
\begin{empheq}
[left={\empheqlbrace}]{align}
\widetilde{\bf x} _t &=F(\widetilde{\bf x}_{t-1}, \widetilde{\bf y} _{t-1}),\label{eq:pred1}\\
\widetilde{\bf y}_t &= h (\widetilde{\bf x} _t), \label{eq:pred2}
\end{empheq}
for $t \in \mathbb{N}$, and initialized with $\widetilde{\bf x}_0 = H^F(\underline{{\bf y}_{-1}})$ and $\widetilde{\bf y}_0 = {\bf y}_0$. If we denote the concatenation of ${\bf y}\in ({\mathbb{R}}^d)^{\mathbb{Z}_{-}} $ and $\mathbf{v} \in ({\mathbb{R}}^d)^{n} $ by ${\bf y} \mathbf{v} = \left(\ldots, {\bf y}_{-1}, {\bf y}_{0}, \mathbf{v}_1, \ldots, \mathbf{v}_n \right) \in ({\mathbb{R}}^d)^{\mathbb{Z}_{-}}  $, we can write our two sequences as
\begin{equation}\label{first steps recursions}
    \widetilde{ {\bf y}} _t=H^F_h( \underline{{\bf y}_{0}}\, \underline{\widetilde{{\bf y}}_{t-1}}), \mbox{ and } \widetilde{{\bf x}}_t = H^F(\underline{{\bf y}_{0}}\, \underline{\widetilde{{\bf y}}_{t-1}}).
\end{equation}
We also let ${\bf x}_t = H^F(\underline{{\bf y}_{t-1}})$. While our aim is to estimate the forecasting error $\Delta {\bf y}_t := {\bf y}_t - \widetilde{\bf y}_t$, we do this by mapping the sequences $\underline{{\bf y}_t}$ and $\underline{ {\bf y} _0} \,\underline{\widetilde{\bf y}_t}$ into the state-space variables ${\bf x}_t, \widetilde{\bf x}_t \in \mathbb{R}^L $ via $H^F$ and studying how the error $\Delta {\bf x}_t := {\bf x}_t - \widetilde{\bf x}_t$ accumulates there. In Theorem~\ref{thm:linearized_seq_bound}, we see that the sequences separate at an exponential rate controlled by the top Lyapunov exponent of the associated dynamical system, and in Theorem~\ref{thm:total_state_error_bound}, we see how this exponential separation in the state space carries over to the predictions of the time series being learnt.

The reason why we called this setup {\it simplified} at the beginning of this paragraph is that, as mentioned also in the discussion of the computability of causal chains, in practice, one does not have at one's disposal a full semi-infinite history of the chain that one wants to forecast but only a finite sample, say $\left\{{\bf y}_{-T_0}, {\bf y}_{-T_0+1}, \ldots, {\bf y}_{0}\right\} $. This fact carries in its wake two complementary sources of error that we briefly describe, together with references where they are specifically addressed: 
\begin{description}
\item [(i)] {\bf The estimation error}: the universality results that we cited in the introduction guarantee the existence of uniform approximants $H _h ^F  $ that satisfy \eqref{condition unif approx} in a variety of situations and for various families of state-space systems. Given a finite sample $\left\{{\bf y}_{-T_0}, {\bf y}_{-T_0+1}, \ldots, {\bf y}_{0}\right\} $, the approximant $H _h ^F  $ in the family is in general estimated using the empirical risk minimization with a prescribed loss. All estimation methods have an error associated that decreases with the sample size. Explicit bounds for this error can be found in \cite{RC10}. 
\item [(ii)] {\bf The state starting error}: as it can be seen in~\eqref{first steps recursions}, the iterations that produce the forecasts are initialized using the state given by $H^F( \underline{{\bf y}_{-1}}) \in {\mathbb{R}}^L$ which is uniquely determined if $F$ has the echo state property and the entire semi-infinite sequence $\underline{{\bf y}_{-1}} $ is known. When only a finite history $\left\{{\bf y}_{-T_0}, {\bf y}_{-T_0+1}, \ldots, {\bf y}_{0}\right\} $ is available, the initial state $H^F( \underline{{\bf y}_{-1}}) $ can be determined only approximately, which leads to what we denominate the {\it state starting error}. If the state map $F$ is contracting in the state variable, then the corresponding filter has what is called the {\it state forgetting property} \cite[Theorem 27]{RC9} elsewhere known as the {\it uniform attraction property} (see \cite{manjunath2022embedding}). This property implies that if we initialize the state equation at an arbitrary state sufficiently far in the past, then the resulting state at $t =0 $ converges to $H^F( \underline{{\bf y}_{-1}}) $ as $T_0 \rightarrow \infty $ regardless of the initialization chosen. This means that in the presence of this property, the state $H^F( \underline{{\bf y}_{-1}}) $ can be approximated by using a sufficiently long washout. In the dynamical systems context  for which an embedding generalized synchronization exists, one can learn a map that allows what is called the {\it cold starting} of the RNN without necessarily using a washout (see \cite{RC26} for a presentation of this methodology). 
\end{description}

\paragraph{Bounds for the forecasting error.} We proceed to develop the necessary theory to estimate the error associated with the iterated multistep forecasting scheme as a function of the forecasting horizon. We do this by mapping the true and predicted sequences into the state space using the state-space system. Predictions in the state space come from the trajectory of the associated dynamical system. Using ergodic theory, we can track the accumulation of error.

Let $L_h$ be a Lipschitz constant for $h$, that is, $\|h({\bf x}^1) - h({\bf x}^2)\| \leq L_h \|{\bf x}^1 - {\bf x}^2\|$ for any ${\bf x}^1, {\bf x}^2 \in \mathbb{R}^L$. According to the identities \eqref{first steps recursions} we have
\begin{multline}
\label{first step ineq z}
\left\|\Delta {\bf y}_t \right\|=\left\|H(\underline{{\bf y}_{t-1}})-H^F_h( \underline{{\bf y}_{0}}\, \underline{\widetilde{{\bf y}}_{t-1}})\right\|\\
=\left\|H(\underline{{\bf y}_{t-1}})-H^F_h(\underline{{\bf y}_{t-1}})+H^F_h(\underline{{\bf y}_{t-1}})-H^F_h( \underline{{\bf y}_{0}}\, \underline{\widetilde{{\bf y}}_{t-1}})\right\|
\leq \varepsilon +\left\|H^F_h(\underline{{\bf y}_{t-1}})-H^F_h( \underline{{\bf y}_{0}}\, \underline{\widetilde{{\bf y}}_{t-1}})\right\|\\
=\varepsilon + \left\|h \left(H^F(\underline{{\bf y}_{t-1}})\right)-h \left(H^F( \underline{{\bf y}_{0}}\, \underline{\widetilde{{\bf y}}_{t-1}})\right)\right\| 
\leq \varepsilon + \left\| h({\bf x}_{t-1}) - h(\widetilde{\bf x}_{t-1})\right\|
\leq \varepsilon + L_h\|\Delta {\bf x}_t\|.
\end{multline}

Thus we can track the prediction error by bounding the error in the state space.
Associated to the state-space system $(F,h)$ is the dynamical system $\Phi\colon \mathbb{R}^L \longrightarrow \mathbb{R}^L$ determined by $\Phi({\bf x}) = F({\bf x}, h({\bf x})).$
According to the iterations \eqref{eq:pred1} and \eqref{eq:pred2}, the sequence $(\widetilde{\bf x}_t)$ is a trajectory of the dynamical system $\Phi$, that is, $\widetilde{\bf x}_t = \Phi(\widetilde{\bf x}_{t-1})$ for $t \in \mathbb{N}$ with initial condition $\widetilde{\bf x}_0 = {\bf x}_0$. Defining ${\bf C}_t := F({\bf x}_{t-1}, {\bf y}_{t-1}) - F({\bf x}_{t-1}, h({\bf x}_{t-1}))$ for $t \in \mathbb{N}$, we have
\begin{align}\label{eq:deltax_error_recursion}
\Delta {\bf x}_t = {\bf C}_{t} + \Phi({\bf x}_{t-1}) - \Phi(\widetilde{\bf x}_{t-1})={\bf C}_{t} + D\Phi(\widetilde{\bf x}_{t-1}) \Delta {\bf x}_{t-1} + {\bf O}_{t},
\end{align}
where ${\bf O}_{t}$ is the error resulting from the Taylor approximation of $\Phi({\bf x}_{t-1}) - \Phi(\widetilde{\bf x}_{t-1})$. Let $L_z$ be a Lipschitz for the second argument of $F$, that is, $\|F({\bf x}, {\bf z}_1) - F({\bf x}, {\bf z}_2)\| \leq L_z \|{\bf z}_1 - {\bf z}_2\|$ for any ${\bf x} \in \mathbb{R}^L$ and any ${\bf z}_1, {\bf z}_2 \in \mathbb{R}^d$. Since $\norm{{\bf y}_{t-1} - h({\bf x}_{t-1})} = \norm{ H(\underline{{\bf y}_{t-2}}) -  H^F_h(\underline{{\bf y}_{t-2}})} \leq \varepsilon$, we can bound $\norm{{\bf C}_t} = \|F({\bf x}_{t-1}, {\bf y}_{t-1}) - F({\bf x}_{t-1}, h({\bf x}_{t-1}))\| \leq L_z \varepsilon$. 

\paragraph{Bounds for the linearized error.} We begin our analysis by ignoring the Taylor approximation error ${\bf O}_{t} $ and considering the linearized sequence
\begin{align*}
{\bf a}_t = {\bf C}_t + D\Phi ( \widetilde{\bf x}_{t-1}) {\bf a}_{t-1}, \quad t \in {\mathbb{N}},\quad {\bf a}_0 = 0.
\end{align*}
Using the ergodic theorems Theorem~\ref{thm:Oseledets_mult_erg_thm} and \ref{thm:barr_pes_Oseledets_extension} introduced in Section~\ref{sec:preliminaries}, we bound the growth of such a sequence in Theorem~\ref{thm:linearized_seq_bound}. The ergodic theorems tell us that, asymptotically, the convergence/divergence of trajectories is dominated by an exponential factor. While this growth estimate may deteriorate over time, it will only do so at a small exponential rate, called the {\it leakage rate} in \cite{berry2022learning}. Including the Taylor error leads to a quadratic term in the sequence of errors. Dealing with the error accumulation in this case is harder, but we bound it up to a finite time horizon later on in Theorem~\ref{thm:total_state_error_bound}. 

\begin{theorem}\label{thm:linearized_seq_bound}
Let $\Phi \in \operatorname{Diff}^1(X)$ be a dynamical system acting on the $q$-dimensional $C^1$ submanifold $X$ of $\mathbb{R}^{L}$ and let $X'$ be a compact invariant subset of $X$. Consider a point ${\bf x} \in X'$ and its trajectory in forward time $(\widetilde{\bf x}_t)_{t \in \mathbb{N}_0}$ under the dynamical system $\Phi$, that is, $\widetilde{\bf x}_t = \Phi(\widetilde{\bf x}_{t-1})$ for $t \in \mathbb{N}$, with the initial condition $\widetilde{\bf x}_0 = {\bf x}$. Along with this, consider a sequence $({\bf C}_t)_{t \in \mathbb{N}}$ such that ${\bf C}_t \in T_{\widetilde{\bf x}_t} X$ for all $t \in \mathbb{N}$ and moreover the sequence is bounded by a constant $K>0$. Construct the sequence ${\bf a} = ({\bf a}_t) \in \left(\mathbb{R}^L\right)^\mathbb{N}$ according to the linear recursion
\begin{align*}
{\bf a}_t = {\bf C}_t + D\Phi ( \widetilde{\bf x}_{t-1}) {\bf a}_{t-1}, \quad t \in {\mathbb{N}},\quad {\bf a}_0 = 0.
\end{align*}
Then for any $\Phi$-invariant Borel probability measure on $X'$, there exists a $\Phi$-invariant set $\Lambda \subset X'$ of full measure such that for every point ${\bf x} \in \Lambda$, for any sufficiently small leakage constant $\delta>0$, the corresponding sequence ${\bf a}$ is bounded by
\begin{equation*}
\|{\bf a}_t\| \leq RK e^{t(\lambda_1^{\text{pos}} + \delta)},
\end{equation*}
where $\lambda_1 = \lambda_1({\bf x})$ is the largest Lyapunov exponent of ${\bf x}$, we write $\lambda_1^{\text{pos}} = \max\{\lambda_1, 0\}$, and $R=R(\delta, {\bf x})$ is a constant depending on the initial point ${\bf x}$ and the leakage constant. 
\end{theorem}

\begin{remark}\label{rem:note on constanst for different values for lyapunov exponents}
\normalfont
We prove in fact that
\begin{equation}\label{ineq:different bounds for different values of Lyapunov exponents}
\|{\bf a}_t\| \leq \begin{cases}
R'R''K e^{t(\lambda_1 + \delta)}, & \mbox{if } \lambda_1-3\delta >0,\\
R'Kt e^{4t\delta}, & \mbox{if } \lambda_1-3\delta =0,\\
R'R''K e^{4t\delta}, & \mbox{if } \lambda_1-3\delta <0,
\end{cases}
\end{equation}
for constants $R' = R'(\delta, {\bf x})$ and $R''=R''(\delta, {\bf x})$. Choosing the constant $R$ appropriately then gives us the bound in the theorem. In particular, for the first and third cases we may take $R(\delta, {\bf x}) := R'(\delta, {\bf x}) R''(\delta, {\bf x})$ and $R(\delta, {\bf x}) := R'(\delta/4, {\bf x}) R''(\delta/4, {\bf x})$, respectively. We note that the second case occurs only when $\lambda_1>0$; in that case, we can ensure by shrinking $\delta$  that we are always in the first case and not in the second one.
\end{remark}

\begin{remark}
\normalfont
In practice, we are concerned with the case when the largest Lyapunov exponent is positive, which is a necessary condition for chaos. In this case, trajectories of the dynamical system that start arbitrarily near to each other may diverge (sensitivity to initial conditions), and we wish to bound the rate at which this occurs. When the leading Lyapunov exponent is negative, trajectories that start near to one another converge, making this scenario trivial. While we deal with all possible values of the leading Lyapunov exponent in our theorem, for nonuniformly hyperbolic behavior to hold, we require that all the Lyapunov exponents be nonzero.
\end{remark}

\begin{remark}
\normalfont
The leakage constant $\delta$ quantifies the rate at which the exponential estimate deteriorates. Typically we take $\delta \ll \min_{i=1,\dots,l({\bf x})} |\lambda_i|$. As noted in the theorem, the constant $R=R(\delta, {\bf x})$ has a dependence on $\delta$. Taking a smaller leakage constant will necessarily entail a larger constant $R$.
\end{remark}

\begin{remark}
\normalfont
We comment on the assumption of the existence of a compact invariant set $X'$. First of all, since the Borel probability measure is defined on $X'$, compactness guarantees that $\Phi$ satisfies the integrability condition in Oseledets' Theorem, the cornerstone of our proof. Secondly, by the Krylov-Bogolyubov theorem, the compactness of $X'$ ensures that a $\Phi$-invariant Borel probability measure for $X'$ exists. Indeed, we can go even further and guarantee the existence of an {\it ergodic} Borel probability measure for $\Phi$, since the set of invariant probability measures is convex in the space of measures on $X'$, and the ergodic measures are precisely the extreme points of this set. For an ergodic invariant probability measure, almost all points will have the same Lyapunov exponent. The question of the uniqueness of such measures is well-studied. Given an ergodic measure $\mu$ on $X'$, there is at most one ergodic invariant measure on $X'$, absolutely continuous with respect to $\mu$. For reference texts on the topic of unique ergodicity see \cite[Chapter 6]{walters2000introduction} and \cite[Chapter 4]{lasota2013chaos}.
\end{remark}

\begin{remark}
\normalfont
The theorem that we just formulated and its proof below corrects a misstep in \cite{berry2022learning}. Indeed, in the proof, we decompose the error accumulation sequence ${\bf a}$ into the Oseledets subspaces and track its growth separately in each space. A problem arises when we need to recombine the errors across the Oseledets subspaces since the angles between Oseledets subspaces may change along a trajectory of the dynamical system so that errors in individual Oseledets subspaces committed at different points along the trajectory do not add together in the same way. We deal with this by employing another property of (Lyapunov) regular trajectories, which bounds the rate at which angles between Oseledets subspaces change.
\end{remark}

\noindent {\bf Proof of the theorem.} 
Let $\Lambda$ be the set of regular points of the dynamical system $\Phi$ in $X'$. As noted in Theorem~\ref{thm:Oseledets_mult_erg_thm}, $\Lambda$ is $\Phi$-invariant and has full measure, for any invariant Borel probability measure. Take ${\bf x} \in \Lambda$. Since $(\widetilde{\bf x}_t)$ is a trajectory of $\Phi$, we have the identity $D\Phi(\widetilde{\bf x}_{t-1}) D\Phi(\widetilde{\bf x}_{t-2}) \cdots D\Phi(\widetilde{\bf x}_s) {\bf C}_s = D\Phi^{t-s}(\widetilde{\bf x}_s)$. Using the initial condition ${\bf a}_0 = 0$ we can write
\begin{align}\label{eq:a_error_expanded}
{\bf a}_t &=  {\bf C}_{t} + \sum_{s=1}^{t-1} D\Phi(\widetilde{\bf x}_{t-1}) D\Phi(\widetilde{\bf x}_{t-2}) \cdots D\Phi(\widetilde{\bf x}_s)\nonumber\\
&=  {\bf C}_{t} + \sum_{s=1}^{t-1} D\Phi^{t-s}(\widetilde{\bf x}_{s}) {\bf C}_{s}
=  {\bf C}_{t} + \sum_{s=1}^{t-1} D\Phi^{t} (\widetilde{\bf x}_0) \left(D\Phi^{s}(\widetilde{\bf x}_{0})\right)^{-1} {\bf C}_{s}
\end{align}
Since ${\bf x}$ is regular, it has Lyapunov exponents $\lambda_1({\bf x})> \cdots > \lambda_{l({\bf x})} ({\bf x})$ and there exists a splitting $T_{\bf x}X = \bigoplus_{i=1}^{l({\bf x})} E_i({\bf x})$ such that
\begin{description}
\item[(i)] $l({\bf x}) = l(\Phi^t({\bf x}))$ for $t \in \mathbb{N}$,
\item[(ii)] $D\Phi^t({\bf x}) E_i({\bf x}) = E_i(\Phi^t({\bf x}))$ for $i = 1, \dots, l({\bf x})$, $t \in \mathbb{N}$, and
\item[(iii)] $\lim_{t \to \infty} \frac{1}{t} \log \norm{D\Phi^t({\bf x}){\bf v}} = \lambda_i({\bf x})$ uniformly over $\{{\bf v} \in E_i({\bf x}) \colon \norm{\bf v} = 1\}.$
\end{description}
Thus by {\bf (ii)} above, for $t \in \mathbb{N}$ we have the splitting $T_{\widetilde{\bf x}_t} X = \bigoplus_{i=1}^{l({\bf x})} D\Phi^t E_i({\bf x})= \bigoplus_{i=1}^{l({\bf x})} E_i(\widetilde{\bf x}_t)$ and we can write ${\bf C}_t = \sum_{i=0}^{l({\bf x})} {\bf C}_t^{(i)}$, where ${\bf C}_t^{(i)} \in E_i(\widetilde{\bf x}_s)$ for $i = 1,\dots,l({\bf x})$. This allows us to define sequences ${\bf a}^{(i)} = ({\bf a}_t^{(i)})$ in the $i$th Oseledets subspace by
\begin{align*}
{\bf a}_t^{(i)} = {\bf C}_t^{(i)} + D\Phi ( \widetilde{\bf x}_{t-1}) {\bf a}_{t-1}^{(i)}, \quad t \in {\mathbb{N}},\quad {\bf a}_0^{(i)} = 0,
\end{align*}
for $i=1,\ldots,l({\bf x})$ so that ${\bf a}_t = \sum_{i=1}^{l({\bf x})} {\bf a}_t^{(i)}$, where ${\bf a}_t^{(i)} \in E_i(\widetilde{\bf x}_t)$ for $i = 1, \dots, l({\bf x})$, $t \in \mathbb{N}_0$. Now fix $i=1,\ldots,l({\bf x})$. By {\bf (iii)} above we have that for any $\delta >0$ there is a $\tau \in \mathbb{N}$ such that for all $t\geq \tau$ and for all ${\bf v} \in E_i({\bf x})$ satisfying $\norm{{\bf v}} = 1,$
\begin{align*}
\lambda_i({\bf x}) - \delta &< \frac{1}{t}\log \norm{D\Phi^t({\bf x}) {\bf v}} < \lambda_i({\bf x}) + \delta,
\end{align*}
and so
\begin{align*} 
e^{t(\lambda_i({\bf x}) - \delta)} &< \norm{D\Phi^t({\bf x}) {\bf v}} < e^{t (\lambda_i({\bf x})+ \delta)}.
\end{align*}
Taking a constant $C(\delta,{\bf x},i)\geq 1$ large enough, we have for all $t \geq 0$, and for any ${\bf v} \in E_i({\bf x})\setminus \{{\bf 0} \},$
\begin{equation*}\label{ineq:C_i}
\frac{1}{C(\delta,{\bf x},i)} e^{t(\lambda_i({\bf x}) - \delta)} \norm{{\bf v}} < \norm{D\Phi^t({\bf x}){\bf v}} < C(\delta,{\bf x},i) e^{t (\lambda_i({\bf x})+ \delta)} \norm{{\bf v}}.
\end{equation*}
Thus, for $t \geq 0$,
\begin{align*}
\norm{\left. D\Phi^t ({\bf x}) \right|_{E_i({\bf x})}} &\leq C(\delta,{\bf x},i) e^{t (\lambda_i({\bf x})+ \delta)},
\end{align*}
and
\begin{align*} 
\norm{\left. \left( D\Phi^t ({\bf x}) \right)^{-1} \right|_{E_i(\widetilde{\bf x}_t)}} &\leq C(\delta,{\bf x},i) e^{-t (\lambda_i({\bf x}) - \delta)}.
\end{align*}
Using the same expansion as in~\eqref{eq:a_error_expanded} for each of the sequences ${\bf a}^{(i)}$, we have
\begin{multline*}
\norm{{\bf a}_t^{(i)}} \leq \norm{{\bf C}_t^{(i)}} + \sum_{s=1}^{t-1} \norm{\left.D\Phi^{t}({\bf x})\right|_{E_i({\bf x})}}\norm{\left.\left(D\Phi^s({\bf x})\right)^{-1}\right|_{E_i(\widetilde{\bf x}_s)}} \norm{{\bf C}_s^{(i)}} \\
\leq \norm{{\bf C}_t^{(i)}} + \sum_{s=1}^{t-1} C(\delta, {\bf x},i)^2 e^{t(\lambda_i({\bf x}) + \delta)} e^{-s(\lambda_i({\bf x}) - \delta)} \norm{{\bf C}_s^{(i)}}
\leq \sum_{s=1}^{t} C(\delta, {\bf x},i)^2 e^{t(\lambda_i({\bf x}) + \delta)} e^{-s(\lambda_i({\bf x}) - \delta)} \norm{{\bf C}_s^{(i)}}.\nonumber
\end{multline*}
Let $C(\delta, {\bf x}) = \max_{i = 1, \dots, l({\bf x})} \{C(\delta, {\bf x}, i) \}$. Then we have
\begin{align}\label{ineq:at_bound_before_angles}
\| {\bf a}_t \| \leq \sum_{i=1}^{l({\bf x})} \| {\bf a}_t^{(i)} \| 
&\leq \sum_{s=1}^{t} \sum_{i=1}^{l({\bf x})} C(\delta, {\bf x},i)^2 e^{t(\lambda_i({\bf x}) + \delta)} e^{-s(\lambda_i({\bf x}) - \delta)} \norm{{\bf C}_s^{(i)}}\nonumber\\
&\leq C(\delta, {\bf x})^2 e^{t(\lambda_1({\bf x}) + \delta)} \sum_{s=1}^{t} e^{-s(\lambda_1({\bf x}) - \delta)} \left( \sum_{i=1}^{l({\bf x})} \norm{{\bf C}_s^{(i)}} \right).
\end{align}
It remains therefore to bound the terms $\sum_{i=1}^{l({\bf x})} \norm{{\bf C}_s^{(i)}}$ by the terms $\| {\bf C}_s \|$. This is an important step which was overlooked in \cite{berry2022learning}. At this point we pause our proof in order to develop the necessary theory to deal with this factor. While in any finite-dimensional space, for linearly independent directions ${\bf v}_1,\ldots, {\bf v}_l$, we can find a constant $C>0$ such that $\sum_{i=1}^l |\alpha_i| \| {\bf v}_i \| \leq C \| \sum_{i=1}^l \alpha_i {\bf v}_i \|$, see for instance \cite[Lemma 2.4-1]{kreyszig1991introductory}, in our case we have to deal with the fact that the angles between the Oseledets subspaces change along the trajectory, so we cannot hope to find a single $C$ which fulfills this identity at every point in time. In order to get around this we need Theorem~\ref{thm:barr_pes_Oseledets_extension}, which tells us how volumes change along trajectories. We begin with a few preliminary observations that allow us to relate volumes and angles between spaces.

\begin{definition}
Let $V, W$ be two subspaces of an inner product space $(H, \langle \cdot,\cdot \rangle)$. We define the angle between $V$ and $W$ as

$$\angle(V, W) = \inf \left\{ \arccos \left( \frac{\abs{\langle {\bf v},{\bf w} \rangle}}{\norm{{\bf v}}\norm{{\bf w}}} \right) \colon {\bf v} \in V, {\bf w} \in W, {\bf v},{\bf w} \neq 0 \right\}.$$
\end{definition}

\begin{remark}\label{rem:sin_cos_sign}
\normalfont
By definition $\angle(V,W) \in [0,\pi/2]$. Thus $\cos \angle(V,W), \sin \angle(V,W) \in [0,1].$
\end{remark}

\begin{remark}\label{rem:angle_cos_spaces}
\normalfont
For any ${\bf v} \in V, {\bf w} \in W$, $\abs{\cos \angle({\bf v},{\bf w})} \leq \cos \angle(V,W).$
\end{remark}
We have the following lemma.

\begin{lemma}\label{lmm:sumnorm_normsum_bound_usingangles}
Let $(H, \langle \cdot,\cdot \rangle)$ be an inner product space. Suppose ${\bf v}_1, \dots, {\bf v}_r \in H$ are linearly independent vectors and define the filtrations $\mathcal{V}\colon V_i:= \text{span} \{ {\bf v}_i, \dots, {\bf v}_r\}, i=1,\dots,r$, and $\mathcal{W}\colon W_i := \text{span}\{ {\bf v}_1, \dots, {\bf v}_i\}$,  $i= 1, \dots, r$. So we have $V_r \subsetneq \cdots \subsetneq V_1$ and $W_1 \subsetneq \cdots \subsetneq W_r$. Further define the angles $\theta_i := \angle(V_{i+1}, W_{i}), i= 1, \dots, r-1$. Then
\begin{align*}
\left\|\sum_{i=1}^r {\bf v}_i\right\| \geq \frac{1}{2^{r-1}} \left( \prod_{i=1}^{r-1} \sin \theta_i \right) \left(\sum_{i=1}^r \norm{{\bf v}_i} \right).
\end{align*}
\end{lemma}

\noindent {\bf Proof of Lemma.} For $r=2$, we take two linearly independent vectors ${\bf v}_1,{\bf v}_2$ forming an angle $\theta$. We denote the angle between the spaces $V_2$ and $W_1$ as defined above by $\theta_1$. By Remark~\ref{rem:angle_cos_spaces} we have
\begin{multline*}
\norm{{\bf v}_1+{\bf v}_2}^2 = \norm{{\bf v}_1}^2 + \norm{{\bf v}_2}^2 + 2 \cos \theta \norm{{\bf v}_1}\norm{{\bf v}_2}    
\geq \norm{{\bf v}_1}^2 + \norm{{\bf v}_2}^2 - \abs{\cos \theta_1} (\norm{{\bf v}_1}^2+\norm{{\bf v}_2}^2)\\
= (1 - \abs{\cos \theta_1}) (\norm{{\bf v}_1}^2+\norm{{\bf v}_2}^2)
\geq \frac{1 - \sqrt{ 1 - \sin^2 \theta_1}}{2} (\norm{{\bf v}_1} + \norm{{\bf v}_2})^2.
\end{multline*}
Now, for $0 \leq a \leq 1$ we have the identity $\frac{1- \sqrt{1- a}}{2} \geq \frac{a}{4}.$
In particular, $\frac{1- \sqrt{1- \sin^2\theta_1}}{2} \geq \frac{\sin^2 \theta_1}{4}$. And so, by Remark~\ref{rem:sin_cos_sign},
\begin{align}\label{ineq:2space_sin}
\norm{{\bf v}_1 + {\bf v}_2} &\geq \frac{\sin \theta_1}{2} (\norm{{\bf v}_1} + \norm{{\bf v}_2}).
\end{align}
We proceed by induction. Suppose the result is true for any set of $r$ linearly independent vectors. Consider a linearly independent set $\{{\bf v}_1, \dots, {\bf v}_{r+1}\}$. Noting that $\left(\Pi_{i=1}^{r-1} \sin{\theta_i}\right) / 2^{r-1} \leq 1$, we use \eqref{ineq:2space_sin} and the induction hypothesis to get
\begin{multline*}
\left\| \sum_{i=1}^{r+1} {\bf v}_i \right\| \geq \frac{\sin{\theta_{r}}}{2} \left( \left\| \sum_{i=1}^{r} {\bf v}_i \right\| + \norm{{\bf v}_{r+1}} \right)
\geq \frac{\sin{\theta_{r}}}{2} \left( \frac{1}{2^{r-1}} \left( \prod_{i=1}^{r-1} \sin \theta_i \right) \left(\sum_{i=1}^{r} \norm{{\bf v}_i}\right) + \norm{{\bf v}_{r+1}} \right)\\
\geq \frac{\sin{\theta_{r}}}{2}  \frac{1}{2^{r-1}} \left( \prod_{i=1}^{r-1} \sin \theta_i \right) \left( \left(\sum_{i=1}^{r} \norm{{\bf v}_i}\right) + \norm{{\bf v}_{r+1}} \right)
= \frac{1}{2^{r}} \left( \prod_{i=1}^{r} \sin \theta_i \right) \left(\sum_{i=1}^{r+1} \norm{{\bf v}_i} \right),
\end{multline*}
as desired. \quad $\blacksquare$

We use Lemma \ref{lmm:sumnorm_normsum_bound_usingangles} to bound $\sum_{i=1}^{l({\bf x})} \norm{{\bf C}_t^{(i)}}$ in terms of $\| {\bf C}_t \|$. Let $\{ {\bf v}_1, \ldots, {\bf v}_N\}$ be any basis for $T_{\bf x}X = \bigoplus_{i=1}^{l({\bf x})} E_i({\bf x})$, such that for each $i=1,\dots,l({\bf x})$, we may pick a basis for $E_i({\bf x})$ from $\{ {\bf v}_1, \dots, {\bf v}_N\}$. Thus this basis is normal with respect to the filtration $\mathcal{V}^+({\bf x})$ defined by \eqref{filtration_forward}. Let the filtrations $\mathcal{V}$ and $\mathcal{W}$ be defined as in Lemma \ref{lmm:sumnorm_normsum_bound_usingangles}, that is, $V_i:= \text{span} \{ {\bf v}_i, \dots, {\bf v}_N\}$, and $W_i := \text{span}\{ {\bf v}_1, \dots, {\bf v}_i\}$, for $i= 1, \dots, N$. Furthermore, for each time step $t \in \mathbb{N}$ define ${\bf v}_i(t) := D\Phi^t ({\bf x}) {\bf v}_i$, for $i=1,\dots,N$. Then we may pick a basis for each $E_i(\widetilde{\bf x}_t)$ from $ \{ {\bf v}_1(t), \dots, {\bf v}_N(t)\} $. Analogous to Lemma \ref{lmm:sumnorm_normsum_bound_usingangles}, we define the filtrations $\mathcal{V}(t) \colon V_i(t) := \text{span} \{ {\bf v}_i(t), \dots, {\bf v}_N(t)\}, i=1,\dots,N$ and $\mathcal{W}(t) \colon W_i = \text{span}\{ {\bf v}_1(t), \dots, {\bf v}_i(t)\}, i=1,\dots,N$, with angles $\theta_i(t) := \angle \left( V_{i+1}(t), W_i(t) \right), i=1,\dots,N-1.$

Now, by Theorem~\ref{thm:barr_pes_Oseledets_extension}, $\lim_{t \to \infty} \frac{1}{t} \log \operatorname{vol}( D\Phi^t({\bf x}) {\bf v}_1,\ldots, D\Phi^t({\bf x}) {\bf v}_{N}) = \sum_{i = 1}^N \lambda^+({\bf x},{\bf v}_i)$. Let $S:= \sum_{i = 1}^N \lambda^+({\bf x},{\bf v}_i)$. Then there exists a constant $C'=C'(\delta, {\bf x} )\geq 1$ such that
\begin{align*}
\frac{1}{C'} e^{t(S-\delta)} < \operatorname{vol}({\bf v}_1(t),& \dots, {\bf v}_N(t) ) < C' e^{t(S+ \delta)}
\end{align*}
Similarly, by definition of the forward Lyapunov exponents, for each $i=1,\dots,N$ we have a constant $C_i=C_i(\delta/N, {\bf x}) \geq 1$ such that for all $t \in {\mathbb{N}},$
\begin{equation*}
\frac{1}{C_i} e^{t(\lambda^+({\bf x},{\bf v}_i)-\delta/N)} < \|{\bf v}_i(t)\|  < C_i e^{t(\lambda^+({\bf x},{\bf v}_i)+ \delta/N)}.
\end{equation*}
Thus 
\begin{align*}
\left(\prod_{i=1}^{N-1} \sin \theta_i(t)\right) \left( \prod_{i=1}^N C_i \right) e^{t(S+\delta)} &> \left(\prod_{i=1}^{N-1} \sin \theta_i(t) \right) \left( \prod_{i=1}^N \norm{{\bf v}_i(t)}\right)\\
&= \operatorname{vol}( {\bf v}_1(t), \dots, {\bf v}_N(t) )>\frac{1}{C'} e^{t(S-\delta)},
\end{align*}
and so
\begin{align*}
\prod_{i=1}^{N-1} \sin \theta_i(t) &> \frac{1}{C' \prod_{i=1}^N C_i} e^{-2t\delta}.
\end{align*}

Now, we may write ${\bf C}_t = \sum_{i=1}^N \alpha_i(t) {\bf v}_i(t)$ for some constants $\alpha_1(t), \dots,\alpha_N(t) \in \mathbb{R}$. By construction, for each $E_i(\widetilde{\bf x}_t)$ we may pick a basis from $ {\bf v}_1(t), \dots, {\bf v}_N(t)$. Thus, writing $\sum_{i=1}^N \alpha_i {\bf v}_i(t) = \sum_{i=1}^{l({\bf x})} {\bf C}_t^{(i)}$, we can group the terms in this sum into their respective Oseledets subspaces. In particular, this implies $\sum_{i=1}^{l(\bf x)} \left\| {\bf C}_t^{(i)} \right\| \leq \sum_{i=1}^N \| \alpha_i(t) {\bf v}_i(t) \|$. By Lemma \ref{lmm:sumnorm_normsum_bound_usingangles}, we have 
\begin{align*}
\left\| {\bf C}_t \right\| &= \left\| \sum_{i=1}^N \alpha_i(t) {\bf v}_i(t) \right\|
\geq \frac{1}{2^{N-1}} \prod_{i=1}^{N-1} \sin \theta_i(t) \left(\sum_{i=1}^N \norm{ \alpha_i(t) {\bf v}_i(t) } \right)\\
&> \frac{1}{2^{N-1} C' \prod_{i=1}^N C_i} e^{-2t\delta} \left(\sum_{i=1}^{l({\bf x})} \norm{{\bf C}_t^{(i)}} \right).
\end{align*}
We obtain the desired bound by returning to \eqref{ineq:at_bound_before_angles}.

\begin{multline*}
\| {\bf a}_t \| \leq C(\delta, {\bf x})^2 e^{t(\lambda_1({\bf x}) + \delta)} \sum_{s=1}^{t} e^{-s(\lambda_1({\bf x}) - \delta)} \left( \sum_{i=1}^{l({\bf x})} \norm{{\bf C}_s^{(i)}} \right)\\
\leq 2^{N-1} C' \left( \prod_{i=1}^N C_i \right) C(\delta, {\bf x})^2 e^{t(\lambda_1({\bf x}) + \delta)} \sum_{s=1}^{t} e^{-s(\lambda_1({\bf x}) - \delta)} e^{2s\delta} \norm{{\bf C}_s}\\
\leq 2^{N-1} C' \left( \prod_{i=1}^N C_i \right) C(\delta, {\bf x})^2 K e^{t(\lambda_1({\bf x}) + \delta)} \sum_{s=1}^{t} e^{-s(\lambda_1({\bf x}) -3 \delta)}.
\end{multline*}
We let 
\begin{align*}
R' &= R'(\delta, {\bf x}) := 2^{N-1} C'(\delta, {\bf x}) \left( \prod_{i=1}^N C_i(\delta/N, {\bf x}) \right) C(\delta, {\bf x})^2, \quad \text{and}\\
R'' &= R''(\delta, {\bf x}):= \frac{e^{3\delta - \lambda_1({\bf x})}}{|1-e^{-\lambda_1({\bf x})+3\delta}|}.
\end{align*}

If $\lambda_1({\bf x}) -3 \delta = 0$, then $\sum_{s=1}^{t} e^{-s(\lambda_1({\bf x}) -3 \delta)} = t$, and $e^{t(\lambda_1({\bf x}) + \delta)} = e^{4t\delta}$, so $\| {\bf a}_t \| \leq R'Kt e^{4t\delta}$. If $\lambda_1({\bf x}) -3 \delta < 0, $ then $\sum_{s=1}^{t} e^{-s(\lambda_1({\bf x}) -3 \delta)} \leq e^{-(t+1)(\lambda_1({\bf x}) -3 \delta)}/(e^{- \lambda_1({\bf x}) +3 \delta}-1)$, so $\| {\bf a}_t \| \leq R'R''Ke^{4t\delta}$. The case we are most concerned with is when $\lambda_1({\bf x}) -3 \delta > 0. $ In this case $\sum_{s=1}^{t} e^{-s(\lambda_1({\bf x}) -3 \delta)} \leq e^{-\lambda_1({\bf x}) +3 \delta}/(1-e^{- \lambda_1({\bf x}) +3 \delta})$, so $\| {\bf a}_t \| \leq R'R''K e^{t(\lambda_1({\bf x}) + \delta)}$. This proves the bounds given in~\eqref{ineq:different bounds for different values of Lyapunov exponents}. Remark~\ref{rem:note on constanst for different values for lyapunov exponents} then completes the proof. \quad $\blacksquare$

\paragraph{The main result.} Theorem~\ref{thm:linearized_seq_bound} deals only with a linearized version of the recursion \eqref{eq:deltax_error_recursion}. We now turn to bounding the actual state error $\Delta {\bf x}_t$. In the bounds above, we simplified the error accumulation by performing a Taylor approximation and ignoring the second-order error incurred. When we factor in the Taylor approximation errors, we can maintain the exponential bound, but only up to a certain finite time threshold. This is made explicit in the following theorem.

\begin{theorem}\label{thm:total_state_error_bound}
Let $H: ({\mathbb{R}}^d)^{\mathbb{Z}_{-}} \longrightarrow {\mathbb{R}}^d $ be a causal chain with infinite memory and let $H^F_h:({\mathbb{R}}^d)^{\mathbb{Z}_{-}} \longrightarrow {\mathbb{R}}^d $ be the functional associated to a state-space system $(F,h) $ that has the echo state property, where $F\colon X \times \mathbb{R}^d \longrightarrow X$ and $h\colon X \longrightarrow \mathbb{R}^d$, for a convex $q$-dimensional $C^2$ submanifold $X$ of $\mathbb{R}^{L}$. Let $\Phi\colon X \longrightarrow X$ be the associated dynamical system, that is, $\Phi({\bf x}) = F({\bf x}, h({\bf x}))$ and suppose $\Phi \in \operatorname{Diff}^2(X)$. Further, let $X'$ be a compact $\Phi$-invariant subset of $X$ and define the set of reachable states of the CCIM $X'':= \{ H^F(\underline{{\bf y}_0}) \mid \underline{{\bf y}_0} \text{ is a realization of the CCIM } $H$\}$. We assume that
\begin{description}
    \item [(i)] there exists an open and convex subset $U$ of $X$ such that $X' \cup X'' \subset U$ and the second order derivatives of $\Phi$ are bounded on $U$, and 
    \item [(ii)] $X'' \cup \Phi(X'')$ is contained in a convex subset of $X$.
\end{description}
Suppose that $F$ is Lipschitz on its second entry with constant $L _z $, uniform with respect to the first entry, that is,\begin{eqnarray*} 
\left\|F(\mathbf{x} , {\bf z}_1)-F(\mathbf{x} , {\bf z}_2)\right\| &\leq L _z \left\|\mathbf{z} _1- \mathbf{z} _2\right\|, \quad \mbox{for all $\mathbf{x} \in X $ and all ${\bf z} _1, {\bf z} _2 \in \mathbb{R}^d  $.}
\end{eqnarray*} 
Assume also that the readout $h: \mathbb{R}^L \longrightarrow {\mathbb{R}}^d $ is Lipschitz with constant $L _h $. Suppose that $H^F_h $ is a uniform $\varepsilon$-approximant of $H$ and hence it satisfies \eqref{condition unif approx}. Let ${\bf y} = ({\bf y}_t) \in  ({\mathbb{R}}^d)^{\mathbb{Z}}$ be a realization of the chain $H$, that is,  ${\bf y}_t= H \left(\underline{{\bf y}_{t-1}}\right) $ for all $t  \in \mathbb{Z} $, and let $\widetilde{ {\bf y}} _t=H^F_h( \underline{{\bf y}_{0}}\, \underline{\widetilde{{\bf y}}_{t-1}}) $, $t \in \mathbb{N} $, be the forecast of the positive terms of ${\bf y} $ using $(F,h) $ and the iterative scheme introduced in~\eqref{eq:pred1} and \eqref{eq:pred2}. Let $\widetilde{\bf x}_t = H^F(\underline{{\bf y}_{0}}\, \underline{\widetilde{{\bf y}}_{t-1}})$ be the corresponding sequence in the state space produced by these recursions and let ${\bf x}_t = H^F(\underline{{\bf y}_{t-1}})$, a sequence in the state space parallel to ${\bf y}$. Then for any $\Phi$-invariant Borel probability measure on $X'$, there is a set $\Lambda \subset X'$ of full measure such that, if ${\bf x}_0 \in \Lambda$, then for sufficiently small $\delta>0$, there exists a constant $R=R(\delta, {\bf x}_0)>0$ and a time horizon $T = T(\delta, {\bf x}_0) \in \mathbb{N}$ such that
\begin{align}
\| {\bf x}_t - \widetilde{\bf x}_t \| &\leq 2\varepsilon L_zR e^{t(\lambda_1^{\text{pos}} + \delta)}\label{ineq:bound_deltax_final}\\
\| {\bf y}_t - \widetilde{\bf y}_t \| &\leq \varepsilon \left( 1 + 2 L_h L_z R e^{t(\lambda_1^{\text{pos}} + \delta)} \right),\label{ineq:bound_deltay_final}
\end{align}
for $t=1,\dots,T$, where $\lambda_1 = \lambda_1({\bf x}_0)$ is the largest Lyapunov exponent of ${\bf x}_0$ under the dynamical system $(X,\Phi)$, and as before we take $\lambda_1^{\text{pos}} = \max\{\lambda_1,0\}.$
\end{theorem}

\begin{remark}
\normalfont
The assumptions {\bf (i)} and {\bf (ii)} above are omitted in \cite{berry2022learning}. However, they are necessary for two reasons: first of all, in the proof below, for a point $\widetilde{\bf x}_t \in X'$ and a left infinite sequence $\underline{{\bf y}_t}$ from the CCIM, we approximate a difference $\Phi(H^F(\underline{{\bf y}_t})) - \Phi(\widetilde{\bf x}_t)$ using Taylor's theorem, which requires both endpoints to be in an open and convex subset of $X$, whence the requirement that $X' \cup X'' \subset U \subset X$ for a convex open set $U$. Requiring that the second order derivatives of $\Phi$ are bounded on $U$ allows us to bound the Taylor approximation from above. Secondly, in the proof we encounter the terms ${\bf C}_t = F({\bf x}_{t-1}, {\bf y}_{t-1}) - F({\bf x}_{t-1}, h({\bf x}_{t-1})) = H^F(\underline{{\bf y}_t}) - \Phi(H^F(\underline{{\bf y}_{t-1}}))$ for a left infinite sequence $\underline{{\bf y}_t}$ from the CCIM. We require that ${\bf C}_t$ be in the tangent space of $X$, which is ensured by the assumption that $X'' \cup \Phi(X'')$ is contained in a convex subset of $X$.
\end{remark}

\begin{remark}
\normalfont
We comment on the embedding requirement in the case when the causal chain consists of observations from a dynamical system. In the paper \cite{berry2022learning}, the map $F$ is assumed to embed the underlying dynamical system in the state-space. This is necessary in order to guarantee the existence of a readout as in Proposition~\ref{GS and conjugacies} that forms a conjugate dynamical system in the state-space. The readout is ultimately the map that is learnt and they work with the corresponding learnt dynamical system as an approximation of the dynamical system conjugate to the original. We circumvent the embedding condition by replacing it with the weaker requirement that the functional $H^F_h$ associated to the state-space system $(F,h)$ uniformly $\varepsilon$-approximates the functional of the causal chain. This is an important improvement as the question of determining when a given state-space system will determine a generalised synchronisation that is an embedding when driven by observations from a dynamical system remains open except for a certain class of linear state-space systems \cite{RC21}. Furthermore, it enables us to extend the results to a much larger class of time series, namely that of CCIMs. CCIMs appear as the deterministic parts of many stochastic processes after using the Wold or Doob decompositions and also as the solutions of infinite dimensional or functional differential equations. While our bounds show an error scaling controlled by the top Lyapunov exponent of the learnt dynamical system, it is interesting to note that many of these applications may not possess an underlying system with Lyapunov exponents at all. 
\end{remark}

\begin{remark}\label{rem:finite_time_horizon}
\normalfont
One of the contributions of this paper is that we deal explicitly and rigorously with the error from the Taylor approximation, in contrast to former efforts, where this error was absorbed by giving the error bounds in terms of a factor $O(e^{n(\lambda_1^{\text{pos}} + \delta)})$. We do this by controlling the error growth within a factor of the exponential term up to a finite time horizon. The expression we find for the forecasting horizon within which we can control the error is
$$T=\left \lfloor \frac{1}{\lambda_1^{\text{pos}} + 5\delta} \log\left( \frac{e^{\lambda_1^{\text{pos}} + 5 \delta} -1}{4\varepsilon L_z L_\Phi K  R e^{8 \delta}} \right) \right\rfloor+1,$$
where $L_\Phi$ is a constant depending on the second order tangent map of $\Phi$ and $K= K(\delta, {\bf x}_0)$ is a constant. When we expand $K(\delta, {\bf x}_0)$ and $R(\delta, {\bf x}_0)$ in terms of the other constants introduced in the proofs we obtain
$$T=\left \lfloor \frac{1}{\lambda_1^{\text{pos}} + 5\delta} \left( \log\left[ \frac{(e^{\lambda_1^{\text{pos}} + 5 \delta}-1)(e^{\lambda_1^{\text{pos}}-3\delta}-1)}{e^{8 \delta}}\right] - \log \left[ {4\varepsilon L_z L_\Phi \left( 2^{N-1} C' \left( \prod_{i=1}^N C_i \right) C^2\right)^2 }\right] \right) \right\rfloor+1.$$ 
It is difficult to make qualitative deductions on how the different variables influence the value of $T$. To begin with, the constants $C, C'$ and $C_i, i=1,\dots,N$ all depend on $\delta$, so while decreasing $\delta$ may serve to increase the first logarithmic term, we don't know how it will influence the second term. Indeed, we expect these constants to increase for a smaller leakage constant, so the influence on $T$ is unknown. A second observation is that all the constants have an implicit dependence on $\varepsilon$, since taking a smaller $\varepsilon$ restricts the set of state-space systems from which we can select $(F,h)$. Therefore, although at first glance decreasing $\varepsilon$ may seem to lead to an increase in the predictable time horizon, this may not necessarily be true in practice.
\end{remark}

\begin{remark}
\normalfont
Another limitation of the bounds in~\eqref{ineq:bound_deltax_final} and \eqref{ineq:bound_deltay_final} is that, while the theorem guarantees the existence of the constant $R=R({\bf x}_0, \delta)$, it does not tell how to calculate it. We propose a numerical procedure for estimating $R$ in Section~\ref{sec:numerial illustration}. 
\end{remark}

\begin{remark}
\normalfont
Notwithstanding what has been said above, in the error bounds \eqref{ineq:bound_deltax_final}-\eqref{ineq:bound_deltay_final}, we do see an interplay between overfitting and approximation error. The top Lyapunov exponent is related to the regularity of the dynamical system generated by the state-space system, while the approximation error $\varepsilon$ determines the set of state-space systems we can choose from. Thus decreasing the approximation error may lead to a highly irregular dynamical system and hence a large leading Lyapunov exponent. See \cite[Remark 12]{berry2022learning} for further discussion on this.
\end{remark}


\begin{remark}\label{rem:stab_gap}
\normalfont
In this setup, it is visible how there are three contributions to the error bound in~\eqref{ineq:bound_deltay_final}, namely the approximation  error $\varepsilon$ of $H$ by $H ^F_h $, the Lipschitz constants $L_h$ and $L_z$, which are related to the structure of the RNN, and the exponential factor $e^{t(\lambda_1^{\text{pos}} + \delta)} $ which is exclusively related to the dynamical features of the map $\Phi \in \operatorname{Diff}^1(X) $ that $(F,h)$ determines. Much work has been done on the ability of state-space systems to reproduce the Lyapunov exponents of the dynamical system being learnt. Two factors come into play here: first, while it is known that under suitable conditions on the GS of a state-space system that embeds the dynamical system in a higher dimensional space, the Lyapunov spectrum of the original system is contained in that of the reconstructed system, a number of additional Lyapunov exponents arise, commonly referred to as {\it spurious Lyapunov exponents}, some of which may exceed the top Lyapunov exponent of the original dynamical system \cite{dechert1996topological,dechert2000largest}; second is the question of continuity of Lyapunov exponents under perturbations of the generating dynamical system, which is discussed in detail in the survey \cite{viana2020dis}. Of particular interest, then, is the {\it stability gap}, namely the difference between the top Lyapunov exponent of the reconstructed system and that of the original system. An upper bound for the stability gap is given in \cite[Theorem 3 (iii)]{berry2022learning}. In our case, however, we are concerned only with the Lyapunov exponents of the dynamical system associated to $(F,h)$. Indeed, it is more useful for us to have bounds for the error in terms of the Lyapunov exponents of the learnt dynamical system, which can be calculated, rather than that of the dynamical system we are learning, which are unknown to us. Nevertheless, work on the stability gap is relevant: as discussed in \cite[Theorem 3]{berry2022learning}, for a state-space system that reconstructs the original dynamical system, the stability gap is always nonnegative. Thus, when learning the observations from a dynamical system, if our trained is conjugate to the underlying dynamical system, we cannot hope to get a reconstructed dynamical system with lower Lyapunov exponent than that of the original. Naturally, due to errors in learning, the trained state-space system will not perfectly reconstruct the dynamical system and we may have a negative stability gap, as illustrated in our numerical section where we train an ESN on the Lorenz system and obtain a maximum Lyapunov exponent lower than that of the true system. In this case, however, we would expect an increase in the other constants involved in the bound, such as $\varepsilon$ and $R$, since if our proxy system has vastly different Lyapunov exponents from that of the true underlying system, it cannot yield a good approximation of it. We note that, in general, CCIMs need not have Lyapunov exponents. 
\end{remark}

\begin{remark}
\normalfont
Let $L_x$ be a Lipschitz constant for the first entry of $F$. It is easy to show that if $F$  and $h$ have bounded differentials and $\lambda _1  $ is the top Lyapunov exponent of the dynamical system $\Phi \in \operatorname{Diff}^1(X) $,  then the Lipschitz constants $L _x,L _h$, and  $L _z$ will satisfy
\begin{equation*}
\label{ineq for Lyapunov}
e^{\lambda _1}\leq L _x+L _h L _z,
\end{equation*}
and we can bound the prediction error as
\begin{equation}
\label{bound iterative ineq}
\left\|\Delta {{\bf y}} _t\right\|\leq  \varepsilon \bigg(1+ L _h L _z\sum_{j=0}^{t-2} \left(L _x+L _h L _z\right)^j\bigg), \quad \mbox{for all $t\in \mathbb{N}$.}
\end{equation}
Now, if the Lipschitz constants can be chosen so that $L _x+L _h L _z < 1 $, then the forecasting error admits a horizon-uniform bound given by the inequality:
\begin{equation*}
\label{bound iterative ineq:uniform}
\left\|\Delta {{\bf y}} _t\right\|\leq \varepsilon \left(\frac{1- L _x}{1- \left(L _x+L _h L _z\right)}\right), \quad \mbox{for all $t\in \mathbb{N}$.}
\end{equation*}
Thus, in the case when $\lambda_1 <0$, using Lipschitz constants may give a better bound than that established using the theory of Lyapunov exponents in~\eqref{ineq:bound_deltax_final} and \eqref{ineq:bound_deltay_final}, which grows exponentially. In the case of a positive maximum Lyapunov exponent, the bound \eqref{bound iterative ineq} grows exponentially with factor $L_x+L_hL_z\geq 1$. Calculations show that in practice this factor is typically much larger than the factor $e^{\lambda_1+\delta}$ using the Lyapunov exponent. Thus ergodic theory provides us with a bound which far outperforms an approach using only tools from analysis.
\end{remark}

\begin{remark}
\normalfont
The Lipschitz constants that appear in the error bounds \eqref{ineq:bound_deltax_final} and \eqref{ineq:bound_deltay_final} are readily computable for most commonly used families of recurrent neural networks. Consider the {\it echo state network} (ESN) \cite{Jaeger04} $(F, h) $ which is defined by using a linear readout map $h(\mathbf{x})= W \mathbf{x}  $, $W \in \mathbb{M}_{d, L} $ and a state map given by 
\begin{equation}
\label{esn reservoir map}
F(\mathbf{x}, {\bf z})=\boldsymbol{\sigma}(A \mathbf{x}+ C {\bf z}+ \boldsymbol{\zeta}), \quad \mbox{with} \quad A \in \mathbb{M}_{L,L}, C \in \mathbb{M}_{L,d}, \boldsymbol{\zeta} \in \mathbb{R}^L, 
\end{equation}  
and where $\boldsymbol{\sigma}: \mathbb{R}^L \longrightarrow [-1,1]^L $ is the map obtained by componentwise application of a squashing function $\sigma: \mathbb{R} \longrightarrow [-1,1] $ that is is non-decreasing and satisfies $\lim_{x \rightarrow -\infty} \sigma(x)=-1 $ and $\lim_{x \rightarrow \infty} \sigma(x)=1 $. Moreover, we assume that $L _\sigma:=\sup_{x \in \mathbb{R}}\{| \sigma' (x)|\}  < +\infty$. With this hypothesis, we note that $\vertiii{D_xF(\mathbf{x}, {\bf z})}\leq L _\sigma\vertiii{A}<+ \infty$ and that $\vertiii{D_zF(\mathbf{x}, {\bf z})}\leq L _\sigma\vertiii{C}<+ \infty$. Consequently, by the mean value theorem we can choose in this case
\begin{equation*}
L _x=L _\sigma\vertiii{A}, \ L_z= L _\sigma\vertiii{C}, \  \mbox{and} \ L _h=\vertiii{W}.
\end{equation*}
We recall that for the ESN family, $L _\sigma\vertiii{A}<1 $ is a sufficient condition for the echo state and the fading memory properties to hold when the inputs are in $\ell^{\infty}(\mathbb{R}^d) $ (see  \cite{RC9}).
With regard to the Lyapunov exponents of the dynamical system associated to the ESN, algorithms exist for computing these.
\end{remark}

\begin{remark}
\normalfont
The approximation bound $\varepsilon> 0 $ in~\eqref{ineq:bound_deltax_final} and \eqref{ineq:bound_deltay_final} which has been chosen so that $\left\|H-H^F_h\right\| _{\infty}< \varepsilon  $ can in practice be put in relation, using Barron-type bounds, with architecture parameters of the RNN families that are being used in the forecasting exercise. Relations of this type can be found in \cite{RC12, RC24}.
\end{remark}

\noindent {\bf Proof of the theorem.} As in the proof of Theorem~\ref{thm:linearized_seq_bound}, we take $\Lambda$ to be the set of (Lyapunov-Perron) regular points in $X'$ of the dynamical system $\Phi$. Let ${\bf x}_0 \in \Lambda$. As noted before, according to the iterations \eqref{eq:pred1}-\eqref{eq:pred2}, the sequence $(\widetilde{\bf x}_t)$ is a trajectory of the dynamical system $\Phi$, that is, $\widetilde{\bf x}_t = \Phi(\widetilde{\bf x}_{t-1})$ for $t \in \mathbb{N}$ with initial condition $\widetilde{\bf x}_0 = {\bf x}_0$. Defining ${\bf C}_t := F({\bf x}_{t-1}, {\bf y}_{t-1}) - F({\bf x}_{t-1}, h({\bf x}_{t-1}))$ for $t \in \mathbb{N}$, we have
\begin{align*}
\Delta {\bf x}_t &= F({\bf x}_{t-1}, {\bf y}_{t-1}) - F({\bf x}_{t-1}, h({\bf x}_{t-1})) + F({\bf x}_{t-1}, h({\bf x}_{t-1})) - F(\widetilde{\bf x}_{t-1}, h(\widetilde{\bf x}_{t-1}))\nonumber\\
&={\bf C}_{t} + \Phi({\bf x}_{t-1}) - \Phi(\widetilde{\bf x}_{t-1})
={\bf C}_{t} + D\Phi(\widetilde{\bf x}_{t-1}) \Delta {\bf x}_{t-1} + {\bf O}_{t},
\end{align*}
where ${\bf O}_{t}$ is the error resulting from the Taylor approximation of $\Phi({\bf x}_{t-1}) - \Phi(\widetilde{\bf x}_{t-1})$. Using the fact that ${\Delta {\bf x}_0} = 0$, we can write 
\begin{multline}\label{eq:deltax_true_error_expanded}
\Delta {\bf x}_t = {\bf C}_{t} + \sum_{s=1}^{t-1} D\Phi(\widetilde{\bf x}_{t-1}) \cdots D\Phi(\widetilde{\bf x}_s) {\bf C}_{s} + \bigg( {\bf O}_{t} + \sum_{s=2}^{t-1} D\Phi(\widetilde{\bf x}_{t-1}) \dots D\Phi(\widetilde{\bf x}_{s}) {\bf O}_s\bigg) \\
=  {\bf C}_{t} + \sum_{s=1}^{t-1} D\Phi^{t-s}(\widetilde{\bf x}_{s}) {\bf C}_{s} + \bigg( {\bf O}_{t} + \sum_{s=2}^{t-1} D\Phi^{t-s}(\widetilde{\bf x}_{s}) {\bf O}_{s}\bigg)
=  {\bf a}_t +  {\bf O}_{t} + \sum_{s=2}^{t-1} D\Phi^{t-s}(\widetilde{\bf x}_{s}) {\bf O}_{s},
\end{multline}
where ${\bf a}$ is the sequence defined by the recursion ${\bf a}_t = {\bf C}_t + D\Phi ( \widetilde{\bf x}_{t-1}) {\bf a}_{t-1}$, for $ t \in {\mathbb{N}}$, with ${\bf a}_0 = 0$, as in Theorem~\ref{thm:linearized_seq_bound}. Since $\norm{{\bf y}_{t-1} - h({\bf x}_{t-1})} = \norm{ H(\underline{{\bf y}_{t-2}}) -  H^F_h(\underline{{\bf y}_{t-2}})} \leq \varepsilon$, using the Lipschitz condition on the second coordinate of $F$, we can bound $\norm{{\bf C}_t} \leq L_z \varepsilon$. Thus there exists a constant $R=R(\delta, {\bf x}_0)$ such that $\| {\bf a}_t \| \leq \varepsilon L_z R e^{t(\lambda_1^{\text{pos}} + \delta)}.$

By Taylor's Theorem we also have $\norm{{\bf O}_t} \leq L_\Phi \norm{\Delta {\bf x}_{t-1}}^2$ where $L_\Phi$ is a constant depending on the second order derivatives of $\Phi$.
Using the same strategies as in the proof of Theorem~\ref{thm:linearized_seq_bound}, we can find a constant $K=K(\delta, {\bf x})\geq 1$ such that
$\norm{D\Phi^{t-s} (\widetilde{\bf x}_s)} \leq K e^{t(\lambda_1^{\text{pos}} + \delta)}e^{-s(\lambda_1^{\text{pos}} - 3\delta)}$. Thus,
\begin{align*}
\norm{\Delta {\bf x}_t} &\leq \varepsilon L_z R e^{t(\lambda_1^{\text{pos}} + \delta)} + \sum_{s=2}^{t-1} K e^{t(\lambda_1^{\text{pos}} + \delta)}e^{-s(\lambda_1^{\text{pos}} - 3\delta)}L_\Phi \norm{\Delta {\bf x}_{s-1}}^2 + L_\Phi \norm{\Delta {\bf x}_{t-1}}^2\\
&\leq e^{t(\lambda_1^{\text{pos}} + \delta)} \bigg( \varepsilon L_z R + L_\Phi K \sum_{s=2}^t e^{-s(\lambda_1^{\text{pos}} - 3\delta)} \norm{\Delta {\bf x}_{s-1}}^2\bigg).
\end{align*}
We can thus construct a sequence ${\bf b} = (b_t)_{t \geq 0}$ of upper bounds for $(\norm{\Delta {\bf x}_t})_{t \geq 0}$ by the recursion
\begin{align*}
b_t = e^{t(\lambda_1^{\text{pos}} + \delta)} \bigg( \varepsilon L_z R + L_\Phi K \sum_{s=2}^t e^{-s(\lambda_1^{\text{pos}} - 3\delta)} b_{s-1}^2\bigg), \quad t \in \mathbb{N}; \quad b_0=0.
\end{align*}
We transform this sequence by dividing out the factor $\varepsilon L_z R \exp{t(\lambda_1^{\text{pos}} + \delta)}$, that is, we form the sequence ${\bf c} = (c_t)$ by $c_t := b_t/ (\varepsilon L_z R \exp{t(\lambda_1^{\text{pos}} + \delta)})$. We have
\begin{align}\label{eq:cseq_corr}
c_t = 1 + \varepsilon L_z L_\Phi K R  e^{-2(\lambda_1^{\text{pos}} + \delta)} \sum_{s=2}^t e^{s(\lambda_1^{\text{pos}} +5 \delta)} c_{s-1}^2 \quad \mbox{for } t \in \mathbb{N}.
\end{align}
We can bound the sequence ${\bf c}$ by an exponential growth factor for a finite time horizon. Consider  \eqref{eq:cseq_corr} and suppose that for $s=1,\dots,t-1$ we have $c_s \leq 1+\eta$, for some $\eta>0$. Then
\begin{align*}
c_{t} &\leq 1 + \varepsilon L_z L_\Phi K R e^{-2 (\lambda_1^{\text{pos}} + \delta)} \sum_{s=2}^t e^{s(\lambda_1^{\text{pos}} +5 \delta)} (\eta+1)^2\\&\leq 1 + (\eta+1)^2 \varepsilon L_z L_\Phi K R \left( \frac{e^{8\delta}}{e^{\lambda_1^{\text{pos}} + 5\delta} -1} \right) e^{(t-1) (\lambda_1^{\text{pos}} + 5\delta)}.
\end{align*}
Thus we can maintain the bound $c_t \leq 1+ \eta$ as long as 
\begin{align*}
t \leq \left \lfloor \frac{1}{\lambda_1^{\text{pos}} + 5\delta} \log\left( \frac{\eta}{(\eta+1)^2} \frac{e^{\lambda_1^{\text{pos}} + 5 \delta} -1}{\varepsilon L_z L_\Phi K R e^{8 \delta}} \right) \right\rfloor+1.
\end{align*}
To maximise the time steps for which this bound holds we optimise $\eta/(\eta+1)^2$. We find

$$\frac{d}{d\eta} \left( \frac{\eta}{(\eta+1)^2}\right) =0$$ for $\eta =1$, so that the function has a maximum of $\frac{1}{4}$. Thus, if $c_s \leq 2$ for $s=1,\dots t-1$, then $c_t \leq 2$, as long as
\begin{align*}
t \leq \left \lfloor \frac{1}{\lambda_1^{\text{pos}} + 5\delta} \log\left( \frac{e^{\lambda_1^{\text{pos}} + 5 \delta} -1}{4\varepsilon L_z L_\Phi K R e^{8 \delta}} \right) \right\rfloor+1.
\end{align*}
that is, we can control the error bound as $\norm{\Delta {\bf x}_t} \leq 2 \varepsilon L_z Re^{t(\lambda_1^{\text{pos}} + \delta)}$ within the aforementioned number of time steps. Finally, by \eqref{first step ineq z}, $\left\| \Delta {\bf y}_t \right\| \leq \varepsilon + L_h\|\Delta {\bf x}_t\| \leq \varepsilon \left( 1 + 2 L_h L_z R e^{t(\lambda_1^{\text{pos}} + \delta)} \right)$. \quad $\blacksquare$

\subsection{The causal embedding forecasting scheme} 
\label{The causal embedding forecasting scheme} 

As an application of the previous results, we compare the {\it causal embedding forecasting scheme} proposed by Manjunath and de Clerq in \cite{manjunath2021universal} with the traditional approach used in reservoir computing and RNNs that we analyzed in the previous section. We start by deriving an upper bound for the prediction error corresponding to Theorem~\ref{thm:total_state_error_bound} for his scheme. 

We start with a brief overview of the causal embedding scheme adapted to our setup. The causal embedding scheme consists of a state-space system that takes the states of a dynamical system as inputs and creates an embedding of the dynamical system in a higher dimensional space. In our setup we only have access to observations from a dynamical system, so we begin by embedding the original dynamical system in a proxy space. For this, we may use a state-space system known to provide embedding, such as Takens' delay embedding scheme. We then apply the causal embedding scheme to the proxy space. Another simplification in our setup is that we assume an invertible dynamical system from the start, a property that is treated only as a special case in \cite{manjunath2021universal}.

Let $\phi \colon M \longrightarrow M$ be a dynamical system, with $M$ a $q$-dimensional compact manifold and $\phi \in \operatorname{Diff}^2(M)$. Let $\omega \colon M \longrightarrow \mathbb{R}^d$ be an observation map and suppose $F \colon \mathbb{R}^{L} \times \mathbb{R}^d \longrightarrow \mathbb{R}^{L}$ is a state-space map whose GS $f \colon M \longrightarrow \mathbb{R}^{L}$ is an embedding, so that the proxy space $M' = f(M)$ is a compact submanifold of $\mathbb{R}^{L}$. By Proposition~\ref{GS and conjugacies}, there exists a readout $h \colon M' \longrightarrow \mathbb{R}^d$ satisfying \eqref{readout forecast} and we can form the conjugate dynamical system $\phi' \colon M' \longrightarrow M'$ given by $\phi'(m') = F(m', h(m'))$ for $ m' \in M'$. We assume that $\phi' \in \operatorname{Diff}^2(M')$.

Whereas in the typical RC setup we would learn the readout $h$, in the causal embedding scheme, we concatenate a second state-space system onto the first. Let $F' \colon \mathbb{R}^{L'} \times M' \longrightarrow \mathbb{R}^{L'}$ be a state-space map with GS $f' \colon M' \longrightarrow \mathbb{R}^{L'}$. We define the map $f_2' \colon M' \longrightarrow \mathbb{R}^{L'} \times \mathbb{R}^{L'}$ by
\begin{equation*}
f_2'(m') = (f'(m'), f'(\phi'(m'))). 
\end{equation*}
Let the image of $M'$ under $f_2'$ be $X$. This is the space of single-delay lags in the state variable of $F'$. We assume an additional property for the state-space map $F'$, namely that of {\it state-input invertibility} (SI-invertibility): $F'$ is SI-invertible when the map $F'_{\bf x}(\cdot) := F'({\bf x}, \cdot) \colon M' \longrightarrow \mathbb{R}^{L'}$ is invertible for all ${\bf x} \in \mathbb{R}^{L'}$. SI-invertibility is not a very strict condition to impose, and can in practice easily be guaranteed. When the map $F'$ is SI-invertible, from a single delay in the state space we can compute the corresponding input from the dynamical system $M'$. This allows us to construct a readout $h' \colon X \longrightarrow M'$ as
\begin{equation*}
h'({\bf x}_{-1}, {\bf x}_0) := \phi' \circ (F'_{{\bf x}_{-1}})^{-1}({\bf x}_0),
\end{equation*}
and we can form a dynamical system $\Phi \colon X \longrightarrow X$ given by $\Phi := f_2' \circ h'$. The map $f_2'$ is an embedding and the dynamical system $(X,\Phi)$ is conjugate to $(M',\phi')$ and hence to $(M,\phi)$. The map $\Phi$ maps a single-delay coordinate one step forward under input from the dynamical system. To see this more clearly, let $U^{F'}\colon (M')^{\mathbb{Z}} \longrightarrow (\mathbb{R}^{L'})^{\mathbb{Z}}$ be the filter corresponding to the state-space map $F'$. Then for a point $m' \in M'$, let ${\bf m}' = (\phi'^t(m'))_{t \in \mathbb{Z}}$ be its trajectory under $\phi'$. For any $t \in \mathbb{Z}$ we have
\begin{align*}
\Phi(U^{F'}({\bf m}')_{t-1}, U^{F'}({\bf m}')_t) &= f_2' \circ \phi' \circ \left( F'_{U^{F'}({\bf m}')_{t-1}} \right)^{-1}(U^{F'}({\bf m}')_{t})\\
&= f_2'\circ \phi' \circ \phi'^{t-1}(m') = \left( f'(\phi'^{t}(m')),f'(\phi'^{t+1}(m')) \right) \\&= ( U^{F'}({\bf m}')_{t}, U^{F'}({\bf m}')_{t+1} ). 
\end{align*}

Essentially, what we have done is to enlarge the state-space map $F'\colon \mathbb{R}^{L'} \times M' \longrightarrow \mathbb{R}^{L'}$ to a map $F'' \colon X \times M' \longrightarrow X$ given by $F''({\bf x}_{-1}, {\bf x}_0, m') = ({\bf x}_0, F'({\bf x}_0, m'))$ for $({\bf x}_{-1}, {\bf x}_0) \in X. $ Thus $\Phi({\bf x}_{-1}, {\bf x}_0) = F''({\bf x}_{-1}, {\bf x}_0, h'({\bf x}_{-1}, {\bf x}_0))$. In the causal embedding scheme, we learn the readout $h'$. Note that we do not require the GS $f'$ to be an embedding but only that the map $F'$ is SI-invertible and has the ESP. The reasoning behind adding a second state-space system onto the first is that the first might be highly irregular, making the readout $h$ hard to learn. By contrast, the state-space map $F'$ may be chosen without requiring its GS to be an embedding, and hence the map $\Phi$, and in particular the readout $h'$, is hoped to be more regular and thus easier to learn.

For a point $m \in M$, let ${\bf y} = ({\bf y}_t)_{t \in \mathbb{Z}}$, where ${\bf y}_t = \omega(\phi^t(m))$, for all $t \in \mathbb{Z}$, be the corresponding sequence of observations. The state space maps $F$ and $F'$ generate sequences ${\bf y}' \in (\mathbb{R}^L)^{\mathbb{Z}}$ and ${\bf x} \in (\mathbb{R}^{L'})^\mathbb{Z}$ according to the iterations
\begin{align}\label{eq:causal_generator}
\left\{ \begin{array}{ll}
            {\bf y}_t' &= F({\bf y}_{t-1}', {\bf y}_{t-1}),\\
            {\bf x}_t &= F'({\bf x}_{t-1}, {\bf y}_{t-1}').
        \end{array}
\right.
\end{align}
We approximate $h'\colon X \longrightarrow M'$ with the learnt map $\widetilde{h}$.
In similar manner to the iterated multistep forecasting scheme discussed previously, we take as input the sequence of past observations $\underline{{\bf y}_{0}} \in (\mathbb{R}^d)^{\mathbb{Z}_-}$ 
so that we only have access to the left infinite sequences $\underline{{\bf y}_0'}$ 
and $\underline{{\bf x}_0}$. 
Initializing $\widetilde{\bf x}_0 = {\bf x}_0$ and $\widetilde{\bf y}_0'= {\bf y}_0'$, for $t \in \mathbb{N}$ we produce predictions $\widetilde{\bf x}_1, \widetilde{\bf x}_2, \dots, \widetilde{\bf x}_t$, $\widetilde{\bf y}_1', \widetilde{\bf y}_2', \dots, \widetilde{\bf y}_t'$, and $\widetilde{\bf y}_1, \widetilde{\bf y}_2, \dots, \widetilde{\bf y}_t$ according to the scheme
\begin{align}\label{eq:causal_prediction}
\left\{ \begin{array}{ll}
            \widetilde{\bf x}_t &= F'(\widetilde{\bf x}_{t-1}, \widetilde{\bf y}_{t-1}'),\\ 
            \widetilde{\bf y}_t' &= \widetilde{h}'(\widetilde{\bf x}_{t-1}, \widetilde{\bf x}_{t}),\\
            \widetilde{\bf y}_t &= f^{-1}(\widetilde{\bf y}_t').
        \end{array}
\right.
\end{align}

\begin{corollary}\label{cor:causal_embedding_error}
Let $\phi \in \operatorname{Diff}^2(M)$ be a dynamical system on a finite-dimensional compact manifold $M$ and $\omega \colon M \longrightarrow \mathbb{R}^d$ an observation map and let the maps $F, f, h, \phi', F', f', f_2', h', \Phi$ and the spaces $M'$ and $X$ be as described above. In particular, the GS $f$ belonging to the state-space map is an embedding so that $M' = f(M)$ is a compact submanifold of $\mathbb{R}^L$, and $\phi' \in \operatorname{Diff}^2(M')$ is conjugate to $\phi$; the state-space map $F'$ is SI-invertible and possesses a GS $f'$, and $\Phi$ is conjugate to $\phi'$ and hence to $\phi$. Suppose that $X$ is a convex $C^2$ submanifold of $\mathbb{R}^{L'}$ and that $\Phi \in \operatorname{Diff}^2(X)$.  Further, let $X'$ be a compact $\Phi$-invariant subset of $X$ and define $X'':= \{ U^{F'}({\bf m}')_0 \colon {\bf m}' = (\phi'^t(m'))_{t \in \mathbb{Z}} \text{ is a trajectory of a point } m' \in M'\}$, the set of reachable states. We assume that
\begin{description}
    \item [(i)] there exists an open and convex subset $U$ of $X$ such that $X' \cup X'' \subset U$ and the second order derivatives of $\Phi$ are bounded on $U$, and that
    \item [(ii)] $X'' \cup \Phi(X'')$ is contained in a convex subset of $X$.
\end{description}
Suppose that $f^{-1}$ is Lipschitz on $M'$ with constant $L_{f^{-1}}$ and $F'$ is Lipschitz on its second entry with constant $L _z' $, uniform with respect to the first entry, that is,
\begin{align*} 
\left\|F'(\mathbf{x} , {\bf z}_1)-F'(\mathbf{x} , {\bf z}_2)\right\| &\leq L _z' \left\|\mathbf{z} _1- \mathbf{z} _2\right\|, \quad \mbox{for all $\mathbf{x} \in X $ and all ${\bf z} _1, {\bf z} _2 \in \mathbb{R}^{N}  $.}
\end{align*} 
Let $\widetilde{h}' \colon X \longrightarrow M'$ be a learnt version of the map $h'$ and let $\varepsilon = \| h' - \widetilde{h}'\|_{\infty}$. Assume that $\widetilde{h}': X \longrightarrow M' $ is Lipschitz with constant $L _{\widetilde{h}'} $. For a point $m \in M$, construct the sequences ${\bf y}, {\bf y}'$, and $ {\bf x}$ as in~\eqref{eq:causal_generator}, and the forecasts $\widetilde{\bf x}_1, \widetilde{\bf x}_2, \dots,\widetilde{\bf x}_t$, $\widetilde{\bf y}_1', \widetilde{\bf y}_2', \dots,\widetilde{\bf y}_t'$, and $\widetilde{\bf y}_1,\widetilde{\bf y}_2, \dots,\widetilde{\bf y}_t$ as in~\eqref{eq:causal_prediction}. Then, for any $\Phi$-invariant Borel probability measure on $X'$, there is a $\Phi$-invariant set $\Lambda \subset X'$ of full measure such that if ${\bf x}_0 \in \Lambda$ the following will hold: For small $\delta>0$, there exists a constant $R=R(\delta, {\bf x}_0)$ and a time horizon $T=T(\delta, {\bf x}_0) \in \mathbb{N}$ such that
\begin{align}
\| {\bf x}_t - \widetilde{\bf x}_t \| &\leq 2\varepsilon L_z' R e^{t(\lambda_1^{\text{pos}} + \delta)}\label{ineq:bound_deltax_causal}\\
\| {\bf y}_t' - \widetilde{\bf y}_t' \| &\leq \varepsilon \left( 1 + 2 L_{\widetilde{h}'} L_z' R e^{t(\lambda_1^{\text{pos}} + \delta)} \right),\label{ineq:bound_deltay'_causal}\\
\| {\bf y}_t - \widetilde{\bf y}_t \| &\leq \varepsilon L_{f^{-1}}\left( 1 + 2 L_{\widetilde{h}'} L_z' R e^{t(\lambda_1^{\text{pos}} + \delta)} \right),\label{ineq:bound_deltay_causal}
\end{align}
for $t=1,\dots,T$, where $\lambda_1 = \lambda_1({\bf x}_0)$ is the largest Lyapunov exponent of ${\bf x}_0$ under the dynamical system $(X,\Phi)$, and as before we take $\lambda_1^{\text{pos}} = \max\{\lambda_1,0\}.$
\end{corollary}

{\bf Proof.}
As noted above, the sequences ${\bf y}'$ and ${\bf x}$ satisfy the state-space equations of $(F'',h')$ where $F'' \colon X \times M' \longrightarrow X$ is given by $F''({\bf x}_{-1}, {\bf x}_0, m') = ({\bf x}_0, F'({\bf x}_0, m'))$ for $({\bf x}_{-1}, {\bf x}_0) \in X. $ Similarly, the sequences $\widetilde{\bf y}_1', \widetilde{\bf y}_2', \dots,\widetilde{\bf y}_t'$, and $\widetilde{\bf x}_1, \widetilde{\bf x}_2, \dots, \widetilde{\bf x}_t$ are generated by the state-space system $(F'',\widetilde{h}')$. Thus, ${\bf y}_{t}' = H^{F''}_{h'} (\underline{{\bf y}_{t-1}}')$ for $t \in \mathbb{Z}$ and $\widetilde{\bf y}_{t}' = H^{F''}_{\widetilde{h}'} (\underline{{\bf y}_0}'\, \underline{\widetilde{\bf y}_{t-1}'})$ for $t \in \mathbb{Z}_+$. Now,
\begin{align*}
\|H^{F''}_{h'} - H^{F''}_{\widetilde{h}'}\|_{\infty} = \|(h' - \widetilde{h}') \circ H^{F''}\|_\infty \leq \|h' - \widetilde{h}'\|_\infty = \varepsilon.
\end{align*}
Thus \eqref{ineq:bound_deltax_causal} and \eqref{ineq:bound_deltay'_causal} follow from Theorem~\ref{thm:total_state_error_bound}. Finally we note that
\begin{align*}
\|{\bf y}_t - \widetilde{\bf y}_t\| = \|f^{-1}({\bf y}_t') - f^{-1}(\widetilde{\bf y}_t')\| \leq L_{f^{-1}} \|{\bf y}_t' - \widetilde{\bf y}_t'\|,
\end{align*}
from which \eqref{ineq:bound_deltay_causal} follows immediately. $\quad \blacksquare$

\section{Numerical illustration}\label{sec:numerial illustration}
\label{Numerical illustration}
In this section we illustrate the error bounds given in Theorem~\ref{thm:total_state_error_bound} for predicting observations coming from dynamical systems using reservoir systems. Our results show the applicability of the error bound to these systems in accurately predicting the rate of growth of the forecasting errors. Nevertheless, there are some limitations to the error bound, which become apparent in our numerical work. First of all, the dynamical systems we learn all have attractors with finite diameter, and hence, while the echo state network (ESN) that we use for the forecasting in \eqref{esn reservoir map} typically reconstructs the attractor, the actual predictions deviate from the true values over time until the error saturates the attractor, being bounded by its diameter, and hence it settles at a constant value. This is a limitation inherent to the system. After saturation, any error bound becomes redundant.

Second is the question which naturally arises from the first limitation, namely: in the context where we deal with observations from an unknown dynamical system, are we able to determine up to what time the error bounds are both valid and meaningful? Theorem~\ref{thm:total_state_error_bound} guarantees validity of the bound up to a finite time horizon $T$. However, as noted in Remark~\ref{rem:finite_time_horizon}, this horizon is difficult to calculate in practice. Notwithstanding this, in the empirical exercises below, it is found that the bound remains valid over the entire forecasted trajectory. More importantly, the fitted bounds in our experiments give us meaningful information about the forecasting error up to saturation. We may calculate the saturation error by running predictions on any two distinct trajectories of the dynamical system with close initial conditions until their separation distance stabilizes. The time horizon up to which error estimates are both valid and useful can then be calculated by checking when the predicted error bound crosses this saturation error.

This brings us to the third though standard limitation, which has to do with the constant $R$ in~\eqref{ineq:bound_deltay_final}, which multiplies the exponential growth term. This factor needs to be calculated in order to fit the error bounds for a given initial condition. The existence of this constant is guaranteed by the use of arguments involving Oseledets' Theorem in the proofs; intuitively, the constant gives the system time to settle until the ergodic properties underlying this theorem take effect. While the theory does not supply an algorithm for calculating $R$, we instead estimate the constant by selecting it as the minimum amongst all possible values of $R$ for which the bound holds. In view of Theorem~\ref{thm:total_state_error_bound}, this may be seen as a minimum estimate for the constant $R$ guaranteed in the theorem. While the constant guaranteed in the theory may be higher, as is often the case in applications, in practice the bound holds for a much lower value of $R$. In the case of reservoir systems that do not learn the given system well, the constant from the theorem may indeed turn out to exceed the error at saturation, and the bound becomes redundant. While on first impression this seems perplexing, this is a practical difficulty inherent to any learning paradigm. Rather than undermining the usefulness of the error bound given in this paper, it serves to demonstrate the inability of the given reservoir computer to track the dynamical system, in which case any error bounds are redundant as saturation is quickly reached.

Despite all these limitations, we emphasize that the usefulness of the bound lies in the fact that the {\it exponential growth rate which ultimately dominates the growth of the error can be calculated}, as it depends on the top Lyapunov exponent of the trained system, to which we have access. Indeed, as the numerical experiments below illustrate, in the case of reservoir systems that demonstrate the ability to produce accurate long-term iterated multistep predictions, the gradient of the log error growth of forecasting errors closely matches the predicted value given in~\eqref{ineq:bound_deltay_final}.

In our numerical illustrations, we consider the Lorenz and R\"{o}ssler systems and train an ESN on an input time series consisting of the full states of the given system. In the case of the Lorenz system, we also consider inputs coming from an observation function on the dynamical system. This forms a causal chain as shown in Corollary \ref{corr:dynamical_systems_and_causal_chains}. The ESNs we use have a state map of the form \eqref{esn reservoir map}, that is, 
\begin{equation}
F(\mathbf{x}, {\bf z})=\boldsymbol{\sigma}(A \mathbf{x}+ \gamma C {\bf z} + \boldsymbol{\zeta}), \quad \mbox{with} \quad A \in \mathbb{M}_{L,L}, \gamma \in \mathbb{R}, C \in \mathbb{M}_{L,d}, \boldsymbol{\zeta} \in \mathbb{R}^L, 
\end{equation}  
with an affine readout given by $h \colon \mathbb{R}^L \longrightarrow \mathbb{R}^d, h({\bf x}) = W {\bf x} + {\bf a}$ for some linear map $W \in \mathbb{M}_{d,L}$ and bias ${\bf a} \in \mathbb{R}^d$. We take the nonlinear function $\boldsymbol{\sigma}(\cdot)$ to be tanh$(\cdot)$ applied componentwise and randomly generate the connectivity matrix $A \in \mathbb{M}_{L, L}$ and input mask $C \in \mathbb{M}_{L, d}$ by selecting entries independently from a uniform distribution over the interval $[-1,1]$, with density $0.01$ and $1$, respectively. The input shift $\boldsymbol{\zeta} \in \mathbb{R}^L$ is simply taken to be zero. Data for trajectories of the different systems is generated by integrating the given dynamical system over a discrete time step $\Delta t$ using the Runge-Kutta45 method. In each case, an input trajectory ${\bf y}_{-T_0}, \dots, {\bf y}_0$ is fed into the randomly generated state map, generating a trajectory in the state space according to the scheme
\begin{equation}
{\bf x}_{t+1} = F({\bf x}_t, {\bf y}_t), \quad \mbox{for}  \quad t = -T_0, \dots, -1, \quad \mbox{where} \quad {\bf x}_{ - T_0} = {\bf 0}.    
\end{equation}
An initial washout of $\tau$ timesteps allows the state space trajectory to converge to the values given by the GS so that the two become indistinguishable, although the state space trajectory generated is arbitrarily initialized at $\bf 0$. This is due to the so-called {\it input forgetting property} \cite{RC9}. The washout data is discarded and a ridge regression with regularization constant $\alpha$ is performed between the input points ${\bf y}_{-T_0+\tau}, \dots, {\bf y}_0$ and state space points ${\bf x}_{-T_0 + \tau}, \dots, {\bf x}_0$, to fit the output parameters $W$ and ${\bf a}$. Predictions are then made by iterating \eqref{eq:pred1} and \eqref{eq:pred2}, that is, 
\begin{empheq}[left={\empheqlbrace}]{align}
\widetilde{\bf x} _t &=F(\widetilde{\bf x}_{t-1}, \widetilde{\bf y} _{t-1}),\\
\widetilde{\bf y}_t &= h (\widetilde{\bf x} _t),
\end{empheq}
for $t = 1, 2, 3, \dots,T_1$ where $\widetilde{\bf x}_0 = {\bf x}_0$ and $\widetilde{\bf y}_0 = {\bf y}_0$. There are several hyperparameters that need to be tuned, namely the regularization constant $\alpha$, the scaling constant $\gamma$ for the input mask, and the spectral radius of the connectivity matrix. Selecting these is a perplexing task as, to our knowledge, no clear analytical method exists by which prediction error as a function of these constants may be minimized. Moreover, the error in prediction is not stable with respect to choices of the constants and, especially in selecting the regularization constant, small changes may lead to largely different results. In order to tune our ESNs individually to each system, we had to resort to a brute force cross validation method minimizing mean squared prediction error over a grid for these parameters. This code is computationally expensive.

\begin{figure}[!htb]
\centering
\includegraphics[scale = 0.35]{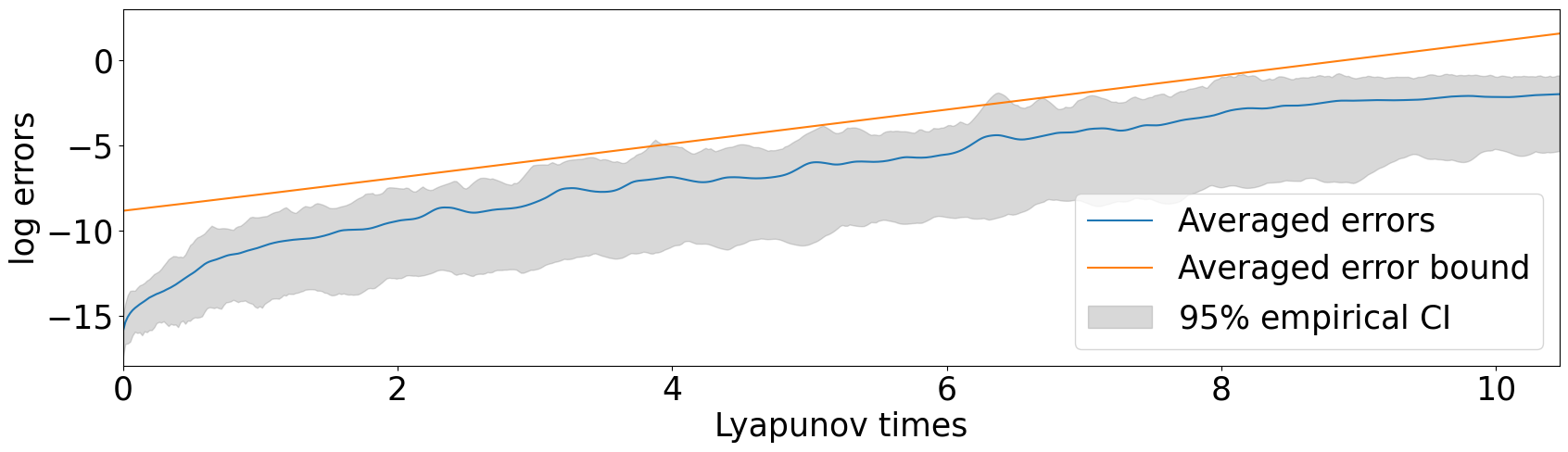}
\caption{Prediction error for the Lorenz system with full access to system states, averaged over 100 trajectories. The error bound plotted has been averaged over the 100 trajectories. The gray band is an emperically calculated $95\%$ confidence interval for the errors committed by the $100$ trajectories. These are all plotted in the log scale.}
\label{fig:lorenz_full_access_errors}
\end{figure}

\begin{figure}[!ht]
\centering
\includegraphics[scale = 0.3]{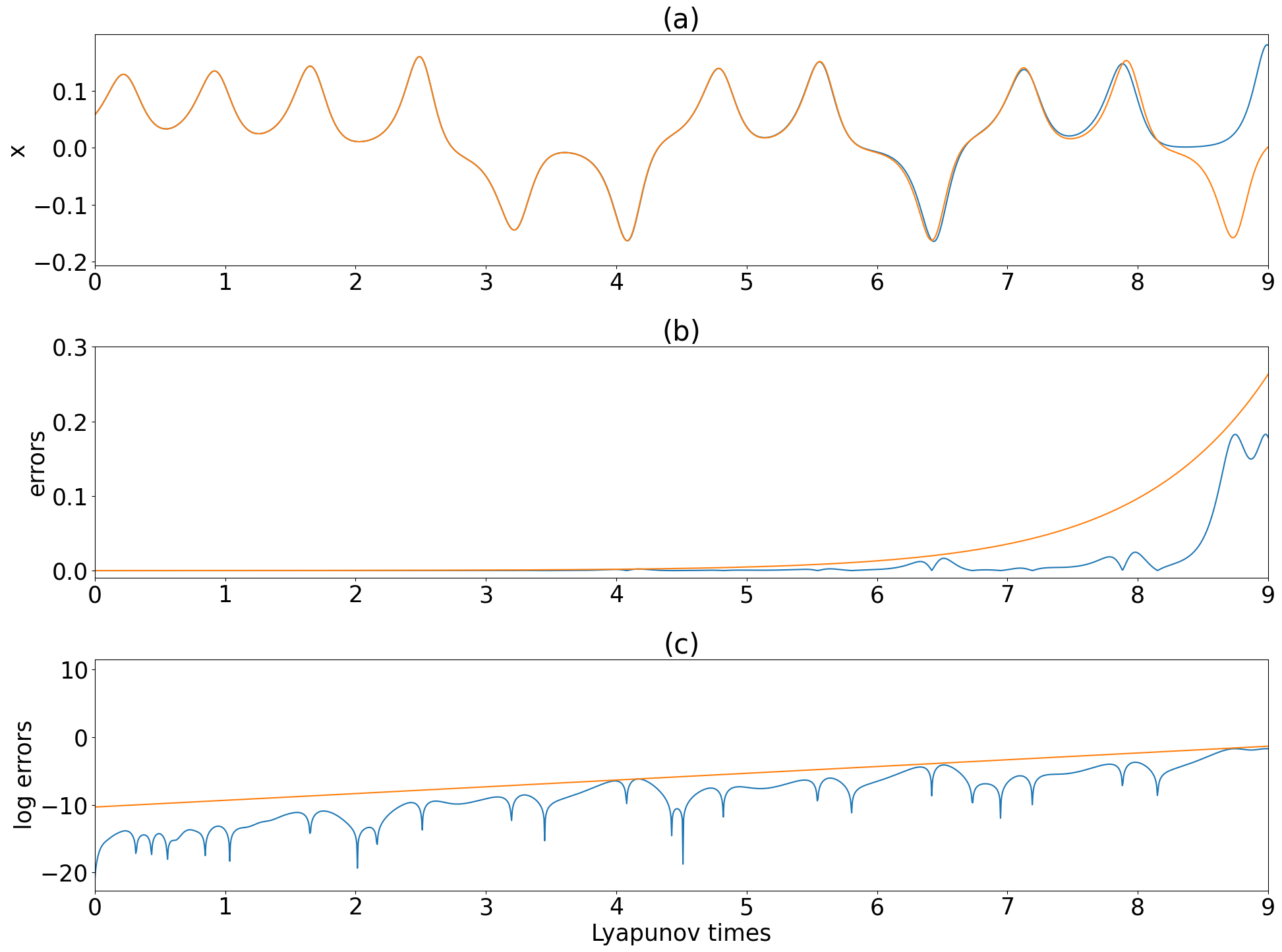}
\caption{Predictions and errors for a single trajectory of an ESN trained on the first coordinate of the Lorenz system.}
\label{fig:lorenz_chain predictions_coordinate}
\end{figure}

\begin{figure}[!ht]
\centering
\includegraphics[scale = 0.35]{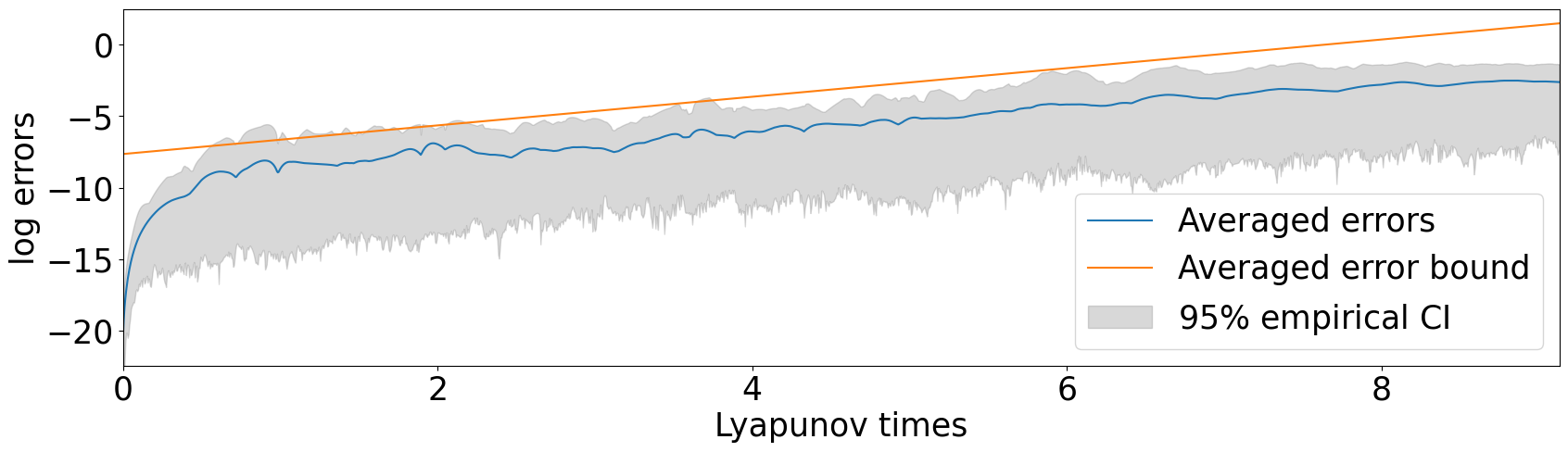}
\caption{Prediction error for the Lorenz system with access to the first coordinate, averaged over 100 trajectories. The error bound plotted has been averaged over the 100 trajectories. The gray band is an emperically calculated $95\%$ confidence interval for the errors committed by the $100$ trajectories. These are all plotted in the log scale.}
\label{fig:lorenz_chain error growth}
\end{figure}

\begin{figure}[!ht]
\centering
\includegraphics[scale = 0.35]{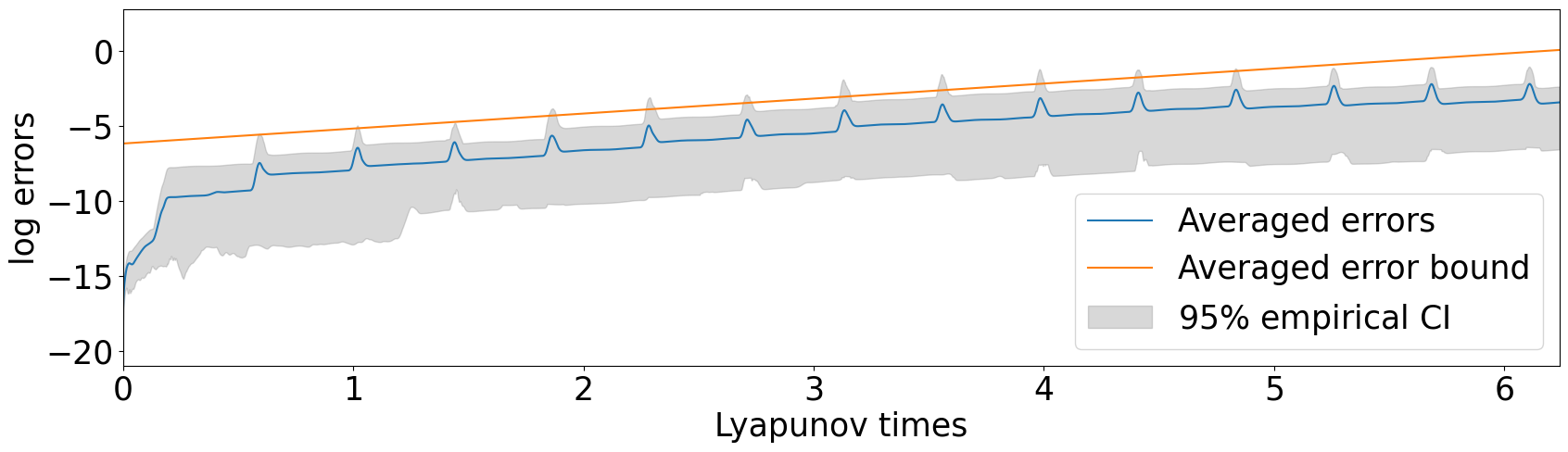}
\caption{Prediction error for the R\"{o}ssler system with full access to system states, averaged over 100 trajectories. The error bound plotted has been averaged over the 100 trajectories. The gray band is an emperically calculated $95\%$ confidence interval for the errors committed by the $100$ trajectories. These are all plotted in the log scale.}
\label{fig:rossler_full access errors}
\end{figure}

\paragraph{Computation of constant $R$.} In our empirical exercises we propose to calculate the constant $R$ for the error bound \eqref{ineq:bound_deltay_final} as follows: for a given trajectory $$({\bf y}_{-T_0}, {\bf y}_{-T_0+1} \dots, {\bf y}_{-1}, {\bf y}_0, {\bf y}_1, {\bf y}_2, \dots, {\bf y}_{T_1}),$$ after the ESN had been trained and predictions $(\widetilde{\bf y}_1, \widetilde{\bf y}_2, \dots, \widetilde{\bf y}_{T_1})$ made, $R$ was selected as

\begin{align}\label{eq:R calculation}
    R := \arg\min_{R'>0} \{R' \colon \varepsilon \left( 1 + 2 L_h L_z R' e^{t \Delta t (\lambda_1^{\text{pos}} + \delta)} \right) \geq \|{\bf y}_t - \widetilde{\bf y}_t \| \quad \text{for all } t =1, \dots, T_1 \}.
\end{align}
\noindent As seen in Figure~\ref{fig:lorenz_chain predictions_coordinate}, this choice of $R$ makes the bound tight, so that the bound touches the actual errors at one point. The value of $\varepsilon$ was estimated from the maximum error in one step ahead predictions of the ESN along the training trajectory $({\bf y}_{-T_0}, {\bf y}_{-T_0 +1}, \dots, {\bf y}_{-1}, {\bf y}_0)$. The Lipschitz constants $L_z$ and $L_h$ can be calculated directly from the trained ESN's parameters, $\delta$ is a constant which we fix by our own choice, and the top Lyapunov exponent was calculated along a single trajectory using Bennetin's algorithm. In the experiments below, $100$ different initial conditions were selected and their trajectories predicted by the ESN. Figures \ref{fig:lorenz_full_access_errors}, \ref{fig:lorenz_chain error growth}, and \ref{fig:rossler_full access errors} show the forecast errors (in blue) for the respective systems, averaged over the $100$ trajectories and plotted in the log scale. The grey band in the figures is a 95$\%$ confidence interval for the forecast errors of trajectories of initial conditions selected from a Gaussian distribution around a point on the respective system's attractor, calculated empirically from the $100$ trajectories at each time step. Note that $R({\bf x}_0, \delta)$ is trajectory dependent, as it depends on the initial condition in the state-space. The value of $R$ used in plotting the error bounds (orange lines) was calculated by taking an average over the values of $R$ given by \eqref{eq:R calculation} for each of the $100$ trajectories. For an individual trajectory, its value of $R$ makes the bound tight, so that it touches the trajectory at at least one point. When we take the average value of $R$ over all the trajectories, this naturally means that some trajectories will, at some points, exceed the plotted error bound. Indeed, this is seen in the figures, where the 95$\%$ confidence interval sometimes crosses the plotted error bound. This does not mean the error bound is invalid, since it holds trajectory-wise, when $R$ is appropriately selected for that trajectory. Positively, the figures demonstrate that selecting $R$ as an average over all the trajectories gives a value for which the error bound is reasonably tight when compared to the average error, while still holding most of the time for 95$\%$ of the trajectories from the selected distribution of initial conditions.

From this observation, we also propose a method for calculating $R$ without making use of the testing data: Typically in the learning task at hand one is given data from the trajectories of the system being learnt, which we divide between trajectories for validation and trajectories for testing. Hyperparameters are selected by a cross validation procedure on the validation data. In this procedure trajectories are divided into a training segment $\{{\bf y}_{-T_0}, {\bf y}_{-T_0+1}, \dots, {\bf y}_0\}$ and a validation segment $\{{\bf y}_1, {\bf y}_2, \dots, {\bf y}_{T_1}\}$. The ESN is initialized with different hyperparameters, its readout trained on the training portion of the trajectories and the ESN then evaluated by its performance in predicting the rest of the trajectories. Finally the best hyperparameter combination is chosen. Following this, the testing trajectories are divided in the same way as the validation trajectories, the second segment this time being for testing. The ESN readout is trained on the training portion of the data and predictions made over the testing segment. To calculate a value of $R$ which would give a reasonable bound for most trajectoires, we may run the ESN with trained readout on the validation data, calculate $R$ for each trajectory using \eqref{eq:R calculation}, and then take an average over these values. As seen in the figures below, with initial conditions coming from the same distribution as those used during validation, the bound accounts for approximately $95\%$ of the trajectories from this distribution.

The main emphasis of the figures remains to illustrate how {\it the gradient of the log error growth matches that given by the bound}, namely the top Lyapunov exponent of the trained system.

\paragraph{The Lorenz system.}
The Lorenz system is given by the three-dimensional first-order differential equation
\begin{empheq}[left={\empheqlbrace}]{align}\label{eq:lorenz system}
\Dot{x} &= -\sigma (x - y),\\
\Dot{y} &= \rho x  - y - xz,\\
\Dot{z} &=  - \beta z + xy,
\end{empheq}
where $\sigma = 10$, $\beta = 8/3 $, and $\rho = 28$.

\begin{figure}[!htb]
\centering
\includegraphics[scale = 0.3]{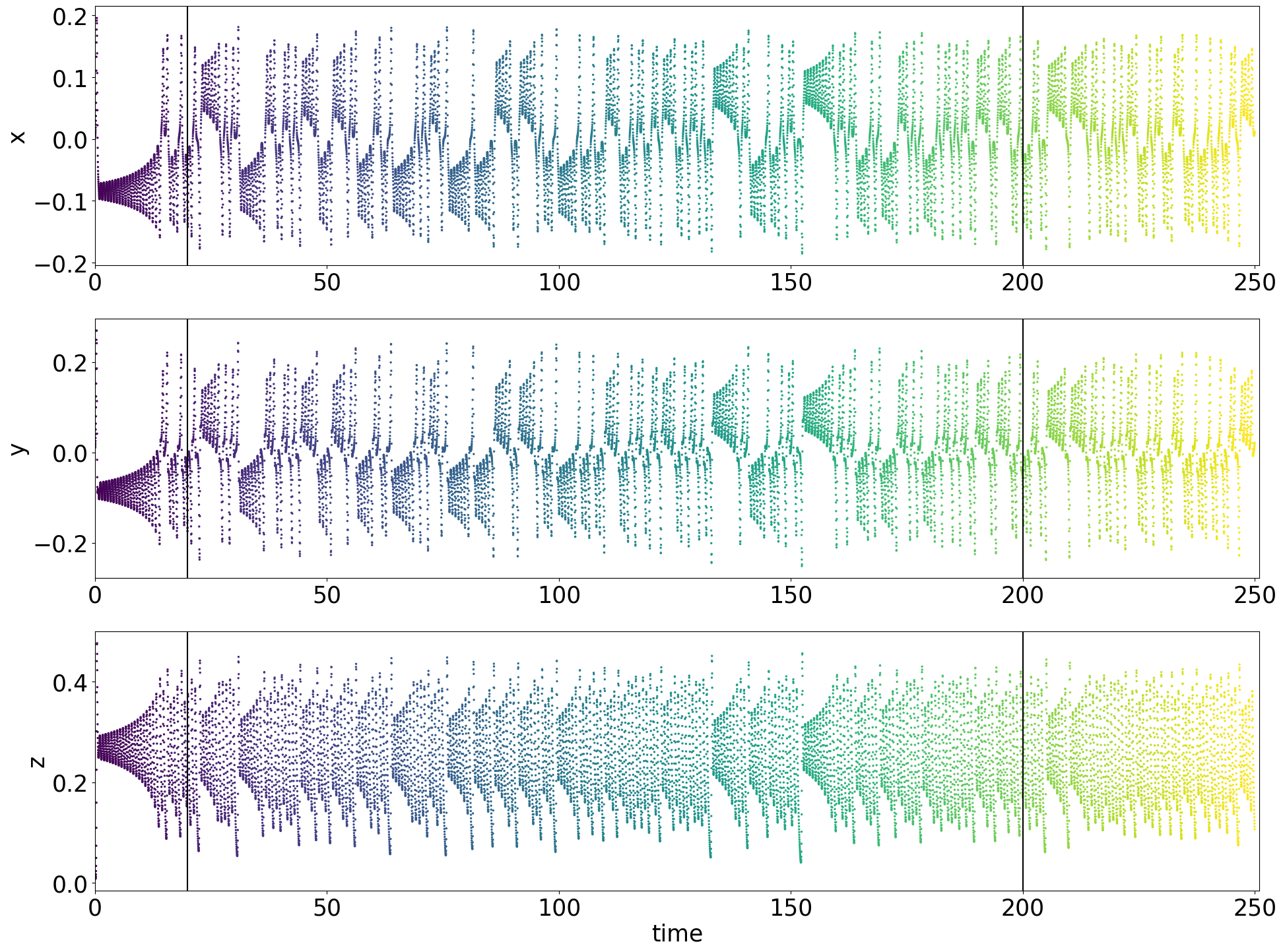}
\caption{Data from a trajectory used for training and forecasting the Lorenz system with access to full coordinatees. The first black line divides the washout from the trajectory data used in training the ESN. The second black line divides the training trajectory from the trajectory used to test forecasting.}
\label{fig:lorenz_full access coordinates}
\end{figure}

To learn the Lorenz system, we chose an ESN with $1000$ neurons whose connectivity matrix had spectral radius $1.2$, and whose input mask was scaled by $\gamma = 13/3$. Inputs were taken as the states of the entire system, sampled at a time step of $\Delta t = 0.02$ over $10000$ time steps. A washout of $1000$ time steps was used. Thereafter, using a Tikhonov regularization, with regularization constant $\alpha = 10^{-14.5}$, between the state space trajectory and desired output trajectory on the remaining $9000$ time steps, the readout $h$ was trained and predictions made for an additional $2500$ time steps. Plotted in Figure~\ref{fig:lorenz_full_access_errors} are the errors in prediction over $10$ Lyapunov times for an RC trained over $100$ trajectories with initial conditions chosen from a Gaussian distribution around $(0,1,1.05)$ with standard deviation $1\times 10^{-10}$, averaged at each point across the $100$ trajectories and plotted in the log scale. The Lyapunov exponent of the RC was calculated at $\lambda_1 = 0.8733$ using Bennetin's algorithm over the length of one trajectory. The value of $\delta$ was fixed at $10^{-3}$ and $\varepsilon = 8.988 \times 10^{-6}$ was calculated by taking the maximum one-step ahead prediction error of the trained ESN over the training portions of the $100$ trajectories. The constant $R$ from \eqref{ineq:bound_deltay_final} was calculated for each trajectory using \eqref{eq:R calculation} and then averaged. Figure~\ref{fig:lorenz_full_access_errors} shows the error bound using this average value, plotted in the log scale. The log errors grow approximately at a linear rate until saturation of the attractor shortly before $9$ Lyapunov times. The gradient of the linear growth is closely matched by the predicted error bound, whose gradient is given by the top Lyapunov exponent of the RC.

\paragraph{The Lorenz system as a causal chain.}
To illustrate our result for a causal chain, we took the observation function $\omega \colon \mathbb{R}^3 \longrightarrow \mathbb{R}, (x,y,z) \mapsto x$, on the Lorenz system. By Corollary \ref{corr:dynamical_systems_and_causal_chains}, sequences of observations through $\omega$ form a causal chain. To learn the Lorenz system using inputs only from the first coordinate, we chose an ESN with $1000$ neurons whose connectivity matrix had spectral radius $1.0$, and whose input mask was scaled by $\gamma = 2.4$. Inputs were sampled at a time step of $\Delta t = 0.005$ over $12000$ time steps. A washout of $1000$ time steps was used. Thereafter, using a Tikhonov regularization, with regularization constant $\alpha = 10^{-14.75}$, between the state space trajectory and desired output trajectory on the remaining $11000$ time steps, the readout $h$ was trained and predictions made for an additional $4000$ time steps. Plotted in Figure~\ref{fig:lorenz_chain error growth} are the errors in prediction over $9$ Lyapunov times for an RC trained on a single trajectory and giving predictions over $100$ trajectories with initial conditions chosen from a Gaussian distribution around $(0,1,1.05)$ with standard deviation $1\times 10^{-10}$, averaged at each point across the $100$ trajectories and plotted in the log scale. The Lyapunov exponent of the RC was calculated at $\lambda_1 = 1.015$ using Bennetin's algorithm over the length of one trajectory. The value of $\delta$ was fixed at $10^{-3}$ and $\varepsilon = 5.733 \times 10^{-7}$ was calculated by taking the maximum one-step ahead prediction error of the trained ESN over the training portions of the $100$ trajectories. The constant $R$ from \eqref{ineq:bound_deltay_final} was calculated for each trajectory using \eqref{eq:R calculation} and then averaged. Figure~\ref{fig:lorenz_chain error growth} shows the error bound using this average value, plotted in the log scale. The log errors grow approximately at a linear rate until saturation of the attractor around $7$ Lyapunov times. The gradient of the linear growth is closely matched by the predicted error bound, whose gradient is given by the top Lyapunov exponent of the RC.

Figure~\ref{fig:lorenz_chain predictions_coordinate} shows the forecasting errors and fitted bound for the $24^{th}$ trajectory from our dataset. The ESN exhibited good learning for this trajectory, giving accurate predictions for about $8$ Lyapunov times. As such it provides an excellent illustration of our error bounds. For this trajectory, the constant $R$ was calculated at $R = 4.218$. Graph (a) shows the true trajectory (in blue) plotted against the forecasted trajectory (in orange). Graph (b) plots the forecast errors (in blue) against the fitted bound (in orange), and Graph (c) is the same, plotted in the log scale. Graph (c) illustrates the similarity between the gradient of the log error bound, given by the top Lyapunov exponent associated to the ESN, and the gradient of the log error growth, demonstrating that, when the ESN has been well-trained, the error bound accurately predicts the rate at which the forecast errors grow.


\paragraph{The R\"{o}ssler system.}

The dynamics of the R\"{o}ssler system is given by

\begin{empheq}[left={\empheqlbrace}]{align}\label{eq:rossler system}
\Dot{x} &= - y - z,\\
\Dot{y} &= x + ay,\\
\Dot{z} &= b + z(x-c),
\end{empheq}
where $a = 0.1, b = 0.1$, and $c = 14.$
To learn the R\"{o}ssler system, we chose an ESN with $1000$ neurons and connectivity matrix with spectral radius $0.2$. The input mask was scaled by $\gamma = 4.65$. Inputs were taken as the states of the entire system, sampled at a time step of $\Delta t = 0.005$ over $23000$ time steps. A washout of $2000$ time steps was used. Thereafter, using a Tikhonov regularization, with regularization constant $\alpha = 10^{-14.5}$, between the state space trajectory and desired output trajectory on the remaining $21000$ time steps, the readout $h$ was trained and predictions made for an additional $30000$ time steps. Plotted in Figure~\ref{fig:rossler_full access errors} are the errors in prediction over $6$ Lyapunov times for an RC trained over $100$ trajectories with initial conditions chosen from a Gaussian distribution around $(2, 1, 5)$ with standard deviation $1\times 10^{-2}$, averaged at each point across the $100$ trajectories and plotted in the log scale. The Lyapunov exponent of the RC was calculated at $\lambda_1 = 0.06933$ using Bennetin's algorithm over the length of one trajectory. The value of $\delta$ was fixed at $10^{-4}$ and $\varepsilon = 9.263 \times 10^{-6}$ was calculated by taking the maximum one-step ahead prediction error of the trained ESN over the training portions of the $100$ trajectories. The constant $R$ from \eqref{ineq:bound_deltay_final} was calculated for each trajectory using \eqref{eq:R calculation} and then averaged. Figure~\ref{fig:rossler_full access errors} shows the error bound using this average value, plotted in the log scale. The log errors grow approximately at a linear rate until saturation of the attractor around $5$ Lyapunov times. The gradient of the linear growth is closely matched by the predicted error bound, whose gradient is given by the top Lyapunov exponent of the RC.

\noindent We summarize a few values from the numerical experiments in Table \ref{tbl:values}.

\begin{table}[!ht]
\centering
\begin{tabular}{|p{3cm}||c|c|c|}
 \hline
 & Lorenz full access & Lorenz $x$ access & R\"{o}ssler full access\\
 \hline
 $\lambda_1$ & $0.8733$ & $1.015$ &$ 0.06933$\\
 $\varepsilon$ & $8.988 \times 10^{-6}$ & $5.733 \times 10^{-7} $ & $ 9.263 \times 10^{-6}$\\ 
 $L_z$ & $82.03$ & $ 44.13$ & $ 88.02$\\
 $L_h$ & $0.2216$ & $ 0.1507$ & $ 0.5116$\\
 $R$ max & $2.597$ & $ 883.7$ & $ 15.46$\\
 $R$ mean & $0.4089$ & $ 62.22$ & $ 2.497$\\
 $R$ sd & $ 0.01407$ & $ 9.700$ & $ 0.03031$\\
 \hline
\end{tabular}
\caption{Values from the numerical experiments.}
\label{tbl:values}
\end{table}

\paragraph{Discussion.}


The numerical experiments illustrate how the error committed using an ESN to predict three different systems grows at an exponential rate, governed by the top Lyapunov exponent associated with the learned ESN. There is an interplay between the Lyapunov exponents of the dynamical system being learnt and those of the trained ESN. This is seen in its effect on the time until saturation. In the case of the Lorenz attractor, which is known to have a top Lyapunov exponent of $0.9057$ \cite{viswanath2000}, the top Lyapunov exponent of the ESN trained on the full system states is closer to that of the true system than the top Lyapunov exponent of the ESN trained on a single coordinate. This is to be expected, given that the former ESN is given access to more information during training. Nevertheless, the ESN trained on a single coordinate does remarkably well, on average reaching saturation error around $7$ Lyapunov times, while the ESN trained on full system states saturates on average at $8$ Lyapunov times.

\noindent Notably, the ESN trained on the full system states has top Lyapunov exponent lower than that of the true Lorenz system. As discussed in Remark~\ref{rem:stab_gap}, decreasing the top Lyapunov exponent of the trained system arbitrarily is expected to lead to an increase in other constants in the bound. The R\"{o}ssler system has a much lower top Lyapunov exponent than the Lorenz system, and hence in learning it is possible to make accurate predictions significantly further into the future. However, when we standardize the time scale in terms of Lyapunov times, the ESNs trained on the Lorenz system seem to outperform that trained on the R\"{o}ssler system, which reaches saturation between $5$ and $6$ Lyapunov times on average. The reasons for this may be a combination of a diverse range of factors from the geometry of the different attractors to the care taken in tuning hyperparameters; at the very least this serves to show that the maximum Lyapunov exponent of a system is not the only factor at play. The performance of a single algorithm may vary across different systems.

\noindent An extensive discussion on the performance of different algorithms in forecasting chaotic systems is beyond the scope of this paper. Notable progress has been made in \cite{gilpin2023modelscaleversusdomain}, where the authors propose a database of $135$ chaotic systems and benchmark standard time-series forecasting methods against larger domain-agnostic algorithms on this database. The same author makes further progress in \cite{zhang2024zeroshotforecastingchaoticsystems}, where algorithms specialized for the chaotic system forecasting task are compared to more general purpose foundation models performing zero-shot learning given limited training data. Another important study is \cite{schotz2024machinelearningpredictingchaotic} where the authors propose a new metric, the `cumulative maximum error', for evaluating the learning performance of different algorithms. They also consider the `valid time' up to which an algorithm makes forecasts whose error remains below a fixed threshold. This is analogous to the time up to saturation in our graphs, putting a well-defined metric on the time up to which an algorithm yields meaningful forecasts.

\noindent Notably, ESNs perform well amongst the algorithms they consider. Since ESNs are able to give predictions that remain valid over a long time scale, this makes them well-suited for our purpose of illustrating the error bounds derived in this paper. While the error bound remains valid even for algorithms that do not exhibit good learning, in such cases the linear growth rate in the log errors is not as visually apparent as saturation is quickly reached. This was seen, for example, in trajectories originating from initial conditions for which the ESNs we trained could not make accurate long term predictions. Nevertheless, averaging over a number of randomly generated initial conditions minimizes the effect of initial conditions for which the predictions were not as good, which tend to be a smaller proportion, and shows that in general the linear growth rate of the log errors is indeed matched by the top Lyapunov exponent of the trained system.

\noindent It was, however, necessary to draw the initial conditions from a distribution concentrated on a small portion of the attractor. The reason for this was that when training the ESN on initial conditions distributed widely across the attractor, learning was not as good and predictions quickly reached saturation. The error bounds given in this paper find their usefulness {\it precisely when a state-space system has been well-trained and is able to make accurate long-term predictions}, since for these systems we may begin to ask further questions, such as at what rate the accuracy of the predictions deteriorates, and up to what time we can trust its predictions. In the case of ill-trained systems, any predictions will soon lose their usefulness. Thus it was necessary to train ESNs that exhibited accurate predictions in the long range.


\section{Conclusions}

In this paper, we have studied the iterated multistep forecasting scheme based on recurrent neural networks (RNN) that are customarily used in reservoir computing for the path continuation of dynamical systems. More specifically, we have formulated and rigorously proved error bounds in terms of the forecasting horizon, functional and dynamical features of the specific RNN used, and the approximation error committed by it. The results that we present show, in particular, the interplay between the universal approximation properties of a family of recurrent networks in the sense of input/output systems and their properties in the realm of the modeling and the forecasting of the observations of autonomous dynamical systems and, more generally, of causal chains with infinite memory (CCIMs). 

The framework in the paper circumvents difficult-to-verify embedding hypotheses that appear in previous references in the literature and applies to new situations like the finite-dimensional observations of functional differential equations or the deterministic parts of stochastic processes to which standard embedding techniques do not necessarily apply. This is an improvement on previous results, which were restricted to sequences consisting of observations of dynamical systems.

The error bounds in this paper show, in particular, that the growth of forecasting error is exponential at a rate governed by the top Lyapunov exponent of the dynamical system associated with the trained RNN in the state space. This is insightful for several reasons: (i) it illuminates the role of state space dynamics in the accuracy of long-term forecasting, (ii) in the case of observations coming from a dynamical system, it illustrates the relationship between the dynamics of the original system and that of the reconstructed dynamical system in the state space and (iii) it provides a tool for comparing the forecasting accuracy of different state-space system learning techniques based on the top Lyapunov exponents of their associated state space dynamical systems in a given learning problem.

As noted in the numerical section, the bounds given have a number of limitations. In particular, the theory does not provide an algorithm for calculating a constant $R$ needed for fitting the bound to a given initial condition, nor does it provide a means for calculating the time horizon up to which the bound remains valid. Instead we propose a numerical means of estimating $R$, and find that, empirically, the bound remains valid over the entire trajectory, providing meaningful information until saturation. In the case of well-trained state space systems used in prediction, the bound demonstrates a striking sharpness, accurately predicting the growth rate of the true forecasting error.

Future work will focus on stochastic extensions of these results in which iterative schemes of this type are used to forecast discrete-time stochastic processes, a context in which RNN universal approximation results have also been formulated.

\bigskip

\noindent {\bf Acknowledgments:} The authors thank G Manjunath for helpful discussions and remarks. JL is funded by a Graduate Scholarship from Nanyang Technological University. JPO acknowledges partial financial support from the School of Physical and Mathematical Sciences of the Nanyang Technological University. LG and JPO thank the hospitality of the Nanyang Technological University and the University of St.~Gallen; it is during respective visits to these two institutions that some of the results in this paper were obtained.

\noindent
\addcontentsline{toc}{section}{Bibliography}
\bibliographystyle{wmaainf}
\bibliography{GOLibrary}
\end{document}